\documentclass[a4paper, 11pt, reqno]{amsart}

\keywords{Tropical geometry, Galois representations, non-Archimedean analytic geometry}
\subjclass[2010]{14T05, 11F23, 14G22}

\usepackage{amsfonts, amsthm, amssymb, amsmath, stmaryrd, mathtools}
\usepackage{mathrsfs,array}
\usepackage{eucal,fullpage,times,color,enumerate,accents}
\usepackage{url}
\usepackage{tikz}
\usetikzlibrary{calc}
\usepackage{xypic}
\usepackage{bm}
\usetikzlibrary{calc}
\usetikzlibrary{fadings}
\usepackage[math]{cellspace}

\definecolor{navy}{rgb}{0,0,.65}

\usepackage[colorlinks,
]{hyperref}
\hypersetup{colorlinks=true,urlcolor=blue,linkcolor=navy,citecolor=navy}

\makeatletter
\def\@tocline#1#2#3#4#5#6#7{\relax
  \ifnum #1>\c@tocdepth 
  \else
    \par \addpenalty\@secpenalty\addvspace{#2}%
    \begingroup \hyphenpenalty\@M
    \@ifempty{#4}{%
      \@tempdima\csname r@tocindent\number#1\endcsname\relax
    }{%
      \@tempdima#4\relax
    }%
    \parindent\z@ \leftskip#3\relax \advance\leftskip\@tempdima\relax
    \rightskip\@pnumwidth plus4em \parfillskip-\@pnumwidth
    #5\leavevmode\hskip-\@tempdima
      \ifcase #1
       \or\or \hskip 1em \or \hskip 2em \else \hskip 3em \fi%
      #6\nobreak\relax
    \dotfill\hbox to\@pnumwidth{\@tocpagenum{#7}}\par
    \nobreak
    \endgroup
  \fi}
\makeatother

\newtheorem{theorem}{Theorem}[subsection]
\newtheorem{corollary}[theorem]{Corollary}
\newtheorem{lemma}[theorem]{Lemma}
\newtheorem{proposition}[theorem]{Proposition}
\newtheorem{definition}[theorem]{Definition}

\newtheorem{questions}{Question}

\newtheorem{construction}[theorem]{Construction}
\newtheorem{example}[theorem]{Example}
\newtheorem{quasi-theorem}[theorem]{Quasi-Theorem}
\newtheorem{blank remark}[theorem]{}
\newtheorem{ssubsection}[theorem]{}

\newtheorem{Th}{Theorem}[]

\newtheorem{rem1}[theorem]{Remark}
\newenvironment{remark}{\begin{rem1}\em}{\end{rem1}}

\newtheorem{not1}[theorem]{Notation}

\newcommand{\CC} {{\mathbb C}}          
            
\newcommand{\GG}{\mathbb{G}}            
		
\newcommand{\PP}{\mathbb{P}}         
\newcommand{\QQ} {{\mathbb Q}}		
\newcommand{\RR} {{\mathbb R}}		
\newcommand{\ZZ} {{\mathbb Z}}		
\newcommand{\TT} {{\mathbb T}}	
\newcommand{\FF}{{\mathbb F}}

\DeclareMathOperator{\Def}{\overset{{}_{\text{def}}}{=}}

\newcommand{\mono}{\!\xymatrix{{}\ar@{^{(}->}[r]&{}}\!}
\newcommand{\epi}{\!\xymatrix{{}\ar@{->>}[r]&{}}\!}

\newcommand{\an}{\text{an}}

\newcommand{\lact}{\raisebox{8pt}{\rotatebox{-90}{$\circlearrowright$}}}

\newcommand{\sm}[1]{\mbox{{\smaller\smaller\smaller\smaller\smaller $ #1 $}}}

\usepackage{graphicx}

\newcommand{\op}{\textsf{{\smaller\smaller\smaller\smaller\smaller op}}}

\begin{document}

\title{Galois actions on analytifications and tropicalizations}

\author{Tyler Foster}
\address{\newline 
L'Institut des Hautes \'Etudes Scientifiques \newline
Le Bois-Marie, 35 route de Chartres \newline
91440 Bures-sur-Yvette \newline
France}
\email{\href{mailto:foster@ihes.fr}{foster@ihes.fr}}

\date{\today}

\pagestyle{plain}

\maketitle

\begin{abstract}
	This paper initiates a research program that seeks to recover algebro-geometric Galois representations from combinatorial data. We study tropicalizations equipped with symmetries coming from the Galois-action present on the lattice of $1$-parameter subgroups inside ambient Galois-twisted toric varieties. Over a Henselian field, the resulting tropicalization maps become Galois-equivariant. We call their images Galois-equivariant tropicalizations, and use them to construct a large supply of Galois representations in the tropical cellular cohomology groups of Itenberg, Katzarkov, Mikhalkin, and Zharkov. We also prove two results which say that under minimal hypotheses on a variety $X_{0}$ over a Henselian field $K_{0}$, Galois-equivariant tropicalizations carry all of the arithmetic structure of $X_{0}$. Namely: (1) The Galois-orbit of any point of $X_{0}$ valued in the separable closure of $K_{0}$ is reproduced faithfully as a Galois-set inside some Galois-equivariant tropicalization of our variety. (2) The Berkovich analytification of $X_{0}$ over the separable closure of $K_{0}$, equipped with its canonical Galois-action, is the inverse limit of all Galois-equivariant tropicalizations of our variety.
\end{abstract}

\setcounter{tocdepth}{1}

\tableofcontents


\begin{section}{Introduction}


\begin{subsection}{Tropicalization: frontiers and confines.}
	Let $K_{0}$ be a non-Archimedean field. Fix an algebraic $r$-dimensional $K_{0}$-torus $T_{0}=\text{Spec}_{\ \!}K_{0}[\ZZ^{r}]$. Each closed subvariety $X_{0}\subset T_{0}$ has an associated polyhedral subspace $\text{Trop}(X_{0})\subset\RR^{r}$ called the {\em tropicalization of $X_{0}$}, and comes equipped with a map on $K_{0}$-rational points $\text{trop}:X_{0}(K_{0})\longrightarrow \text{Trop}(X_{0})$ that takes each $K_{0}$-rational $r$-tuple $x=(x_{1},\dots,x_{r})$ in $X_{0}$ to its coordinatewise valuation
	$$
	\text{trop}(x)\Def
	\big(-\log|x_{1}|,\dots,-\log|x_{r}|\big)
	\ \in\ 
	\RR^{r}.
	$$
The map extends to a surjective map on the Berkovich analytification
	$$
	\text{trop}:X_{0}^{\an}\epi\text{Trop}(X_{0}).
	$$
In this way, tropicalization is an operation that collapses large swaths of information in the algebraic variety $X_{0}$, reducing the variety to a far more tractable object. The beauty of tropicalization lies in the fact that despite eliminating so much information, it preserves important algebro-geometric features of $X_{0}$. The guiding example of this sort is G. Mikhalkin's enumerative result \cite[\S3, Theorem 1]{Mik:03} \cite[\S7.1, Theorem 1]{Mik:05}. Motivated by ideas of M. Kontsevich \cite{Kontsevich:00}, Mikhalkin showed that Gromov-Witten invariants in $\PP^{2}_{\sm{\CC}}$ can be recovered as counts of tropical curves. In a similar vein, B. Osserman, S. Payne, and J. Rabinoff showed that intersection numbers between algebraic varieties can be recovered via intersection numbers of their tropicalizations \cite{Rab:12} \cite{OP:13} \cite{OR:13}.
	
	Given this situation, one is inclined to ask if there are other important geometric features of algebraic varieties that are preserved under tropicalization. Having fixed a separable closure $K\Def K^{\text{sep}}_{0}$, one obvious candidate is the canonical action of the absolute Galois group
	$$
	G\ \Def\ \text{Gal}({K}/K_{0})
	$$
on the variety
	$$
	{X}
	\ \Def\ 
	X_{0,{K}}.
	$$
Unfortunately, there is an immediate obstacle to the recovery of Galois actions in tropicalizations. Indeed, if $K_{0}$ is Henselian, then there exists a unique absolute value on ${K}$ extending the absolute value on $K_{0}$ \cite[Lemma 4.1.1]{EP:05}. Given a finite Galois extension $L/K_{0}$ with Galois group $H\Def\text{Gal}(L/K_{0})$, the restriction of this absolute value to $L$ is given by
	$$
	|a|_{L}
	\ \Def\ 
	\Big|\prod_{h\in H}h(a)\Big|^{1/d},
	$$
where $d=\#H$. From the form of $|-|_{L}$ it follows that the $G$-action on $L$-rational points of $X$ becomes trivial under the tropicalization map
	\begin{equation}\label{equation: coordinatewise tropicalization for separable closure}
	\text{trop}:{X}({K})\longrightarrow\text{Trop}({X}).
	\end{equation}
As we range over all Galois extensions of $K_{0}$, we see that \eqref{equation: coordinatewise tropicalization for separable closure} collapses all Galois orbits, and one can show that the same is true for the Berkovich analytification of ${X}$; all Galois orbits collapse under tropicalization.
	
	This negative observation has lead to a pervasive idea that Galois actions on algebraic varieties are forever lost to tropical geometers. But this is not the case. The tool needed for tropicalizing algebraic varieties in a manner that is sensitive to their Galois structures is already present in the classification of non-split tori \cite[Chapter III, \S8.12]{Borel}, which plays a central role in the Borel-Chevalley theory of linear algebraic groups, completed in the late `50s and early `60s  \cite{Ch:05} \cite{BT:65} \cite{Bor:69}.
\end{subsection}
	

\begin{subsection}{Galois actions on tropicalizations}\label{subsection: Galois actions on tropicalizations}
	Recall that a (not necessarily split) {\em algebraic torus} over $K_{0}$ is any algebraic $K_{0}$-group $T_{0}$ such that for some finite separable extension $L/K_{0}$, the algebraic $L$-group $T_{0,L}$ is isomorphic to the split algebraic torus $\GG^{r}_{\text{m},L}=\text{Spec}_{\ \!}L[\ZZ^{r}]$ for some $r$. Since $T\Def T_{0,K}$ comes with a canonical $G$-action, its character lattice $M=M(T)\Def\text{Hom}_{\sm{{K}\text{-}\bold{Grp}}}({T},\mathbb{G}_{\text{m},\sm{{K}}})$ comes equipped with a continuous $G$-action
	$$
	G
	\ \lact\ 
	M.
	$$
This gives rise to a functor
	\begin{equation}\label{equation: functor on free lattices}
	\!\!\!\!
	\scalebox{.9}{$
	\left\{\!\!
	\begin{array}{rl}
	\mbox{{\bf objects:}} & \!\!\!\mbox{algebraic $K_{0}$-tori}
	\\
	\mbox{{\bf morphisms:}} & \!\!\!\mbox{algebraic $K_{0}$-group}
	\\
	& \!\!\!\mbox{homomorphisms}
	\end{array}
	\!\!\right\}
	\ 
	\xrightarrow{\ \ T_{0}\ \mapsto\ M(T)\ \ }
	\ 
	\left\{\!\!
	\begin{array}{rl}
	\mbox{{\bf objects:}} & \!\!\!\mbox{finite rank lattices with}
	\\
	& \!\!\!\mbox{continuous $G$-action}
	\\
	\mbox{{\bf morphisms:}} & \!\!\!\mbox{$G$-equivariant}
	\\
	& \!\!\mbox{homomorphisms}
	\end{array}
	\!\!\right\}.
	$}
	\end{equation}
The fundamental classification result for algebraic tori over $K_{0}$ says that the functor \eqref{equation: functor on free lattices} is an equivalence of categories \cite[\S III.8.12]{Borel}. The quasi-inverse $T_{0}(-)$ takes a finite rank lattice $M$ with Galois-action $G\ \lact\ M$, forms the split torus ${T}(M)\Def\text{Spec}_{\ \!}{K}[M]$ with twisted $G^{\op}$-action\footnote{Throughout the paper, {\em $G$-action} is synonymous with {\em left $G$-action}, and {\em $G^{\op}$-action} is synonymous with {\em right $G$-action}. Maps equivariant with respect to either $G$- or $G^{\op}$-actions are called {$G$-equivariant}. Given $g\in G$, we write ``$g(-)$" to denote its application with respect to a $G$-action, and ``$g^{\ast}(-)$" to denote its application with respect to a $G^{\op}$-action. The general rule of thumb throughout the paper is: $G$-actions occur on spaces of functions, $G^{\op}$-actions occur on underlying topological spaces.} coming from the Galois-action on $M$, and then produces an algebraic torus $T_{0}(M)$ over $K_{0}$ via Galois descent:
	$$
	T_{0}(M)
	\ \ \ \Def\ \ \ 
	T(M)
	\big/
	G.
	$$
	
	If $X_{0}$ is a $K_{0}$-variety, then we can ask for closed embeddings of $X_{0}$ into (not necessarily split) algebraic tori over $K_{0}$, or what is the same, $G$-equivariant closed embeddings of ${X}=X_{0,{K}}$ into $G^{\op}$-twisted tori ${T}(M)=\text{Spec}_{\ \!}{K}[M]$. Given such an embedding
	$
	\imath:{X}\mono {T}(M),
	$
the tropicalization of ${X}$ with respect to $\imath$ is a polyhedral subspace that we denote
	$$
	\text{Trop}({X},\imath)\ \subset\ N_{\RR}\Def\text{Hom}_{\ZZ}(M,\RR).
	$$
The $G$-action on $M$ puts extra structure on $\text{Trop}({X},\imath)$. First, it induces a $G^{\op}$-action on $N_{\RR}$.
Recall that the tropicalization map
	\begin{equation}\label{equation: trop to N}
	\text{trop}:{X}^{\an}\ \Def\ (X_{0,K})^{\an}\longrightarrow N_{\RR}
	\end{equation}
takes each multiplicative seminorm $|-|_{x}:{K}[M]\longrightarrow\RR_{\ge0}$, constituting a point of ${X}^{\an}$, to the composite
	$$
	M\mono{K}[M]\xrightarrow{\ \ |-|_{x}\ \ }\RR_{\ge0}\xrightarrow{\ -\log\ }\RR.
	$$
When $K_{0}$ is Henselian, the fact that $\imath$ is $G$-equivariant implies that the subspace $\text{Trop}({X},\imath)\subset N_{\RR}$ is $G^{\op}$-invariant and that the map \eqref{equation: trop to N} is $G$-equivariant. For details, see Lemmas \ref{lemma: invariance and equivariance} and \ref{lemma: invariance and equivariance}. To highlight the fact that $\text{Trop}(X,\imath)$ comes equipped with a $G^{\op}$-action, we denote it
	$$
	\text{Trop}_{G}({X},\imath),
	$$
and refer to it as a {\em Galois-equivariant} or {\em $G$-equivariant tropicalization of $X$}. We refer to the resulting $G$-equivariant map
	$$
	\text{trop}:{X}^{\an}\epi\text{Trop}_{G}({X},\imath)
	$$
as the corresponding {\em $G$-equivariant tropicalization map}.

\begin{remark}
	See Examples \ref{example: Brauer-Severi B}, \ref{Galois equi trop of P^1}, and \ref{example: Brauer-Severi C} below for pictures of Galois-equivariant tropicalizations.
	
	Although described above only for Galois-invariant closed embeddings of $X$ into Galois-twisted tori, our construction of Galois-equivariant tropicalizations extends to all Galois-equivariant closed embeddings of $X$ into Galois-twisted toric varieties. One way these closed embeddings arise naturally is from closed subvarieties of arithmetic toric $K_{0}$-varieties in the sense of E. J. Elizondo, P. Lima-Filho, F. Sottile, and Z. Teitler \cite[\S3]{ELFST}. We develop the extended theory of Galois-equivariant tropicalization in \S\ref{section: Galois equivariant tropicalizaiton}.
\end{remark}

\end{subsection}
	

\begin{subsection}{First questions: faithfulness.}\label{subsection: first questions}
	Given the existence of Galois-equivariant tropicalizations, one would like to understand just how faithfully they reproduce the Galois-action on a given algebraic variety. Because tropicalization collapses so much of the information inside an algebraic variety, there are various concrete questions one can ask.

\begin{questions}\label{questions: recovering orbits tropically}
\normalfont
	Assume that $K_{0}$ is Henselian. Given a point $x\in{X}^{\an}$ with finite $G^{\op}$-orbit, for instance a point $x\in X({K})$, does there exist a $G$-equivariant tropicalization $\text{Trop}_{G}({X},\imath)$ such that the $G^{\op}$-set $G^{\op}x\subset X^{\an}$ maps bijectively onto its image under the $G$-equivariant composite map
	$$
	\ \ 
	G^{\op}x
	\mono
	{X}^{\an}
	\epi
	\text{Trop}_{G}({X},\imath)
	\ \ ?
	$$
\end{questions}
	
\begin{remark}\label{remark: other geometric structures}
	We provide an affirmative answer to Question \ref{questions: recovering orbits tropically} in Theorem \ref{theorem: main C} below under minimal hypotheses on $X_{0}$ when $K_{0}$ is a perfect field, or when $X$ is quasiprojective. This suggests that one should ask Question \ref{questions: recovering orbits tropically} with ``point" replaced by other geometric structures on ${X}$. We plan to address variants of Question \ref{questions: recovering orbits tropically} in a future paper. One concrete variant appears in Question \ref{questions: recovering representations on cohomology} below. Let us first offer a global version of Question \ref{questions: recovering orbits tropically}.
\end{remark}

\begin{questions}\label{questions: recovering Berkovich analytification tropically}
\normalfont
	Can we reconstruct the Berkovich analytification ${X}^{\an}$ {\em together with its natural $G^{\op}$-action} from Galois-equivariant tropicalizations of ${X}$?
\end{questions}
	
\begin{remark}
	We provide an affirmative answer to Question \ref{questions: recovering Berkovich analytification tropically} in Theorem \ref{theorem: main B} below, again under minimal hypotheses on $X_{0}$ when $K_{0}$ is a perfect field, or when $X$ is quasiprojective.
\end{remark}

\begin{questions}\label{questions: recovering representations on cohomology}
\normalfont
	Fix a prime $\ell$. Can we realize the $\ell$-adic representations $G\ \lact\ H^{n}({X}_{\text{\'et}},\QQ_{\ell})$ using representations of $G$ on appropriately defined tropical cohomology groups $H^{n}_{\text{trop}}\big(\text{Trop}({X},\imath),\QQ_{\ell}\big)$?
\end{questions}
	
\begin{remark}
	We do not know if Question \ref{questions: recovering representations on cohomology} has an affirmative answer. However, we show in Theorem \ref{theorem: main D} and Remark \ref{remark: large supply of G-representations} below that the tools we develop to answer Questions \ref{questions: recovering orbits tropically} and \ref{questions: recovering Berkovich analytification tropically} allow us to produce a large supply of representations of $G$ in the tropical cellular cohomology groups of I. Itenberg, L. Katzarkov, G. Mikhalkin, and I. Zharkov \cite{IKMZ:16}. We plan to address Question \ref{questions: recovering representations on cohomology} further in a future paper.
\end{remark}

\end{subsection}
	

\begin{subsection}{Main results}\label{section: main results}
	Fix a $K_{0}$-variety $X_{0}$. Let $X=X_{0,{K}}$ with its canonical $G^{\op}$-action over $\text{Spec}_{\ \!}K_{0}$. Let $\bold{TEmb}_{G}({X})$ denote the category where an object is any $G$-equivariant closed embedding into a $G^{\op}$-twisted toric variety
	$
	\imath:{X}\mono Y_{\Sigma}
	$
(see Definition \ref{definition: twisted toric variety}), and where a morphism from one $G$-equivariant closed embedding $\imath:{X}\mono Y_{\Sigma}$ to another $\imath':{X}\mono Y_{\Sigma'}$ is any morphism of $G^{\op}$-twisted toric varieties $f_{\varphi}:Y_{\Sigma}\longrightarrow Y_{\Sigma'}$ (see Definition \ref{definition: twisted toric variety}) fitting into a commutative diagram
	$$
	\begin{xy}
	(0,0)*+{\ {X}\ }="1";
	(17,8)*+{Y_{\Sigma}}="2";
	(17,-8)*+{Y_{\Sigma'}.\!}="3";
	{\ar@{^{(}->}^{\imath} "1"; "2"};
	{\ar@{_{(}->}_{\imath'} "1"; "3"};
	{\ar^{f_{\varphi}} "2"; "3"};
	\end{xy}
	$$

\begin{definition}\label{defintion: G-equivariant system}
\normalfont
	A {\em system of $G$-equivariant toric embeddings of ${X}$}, denoted $\mathcal{S}_{G}$, is any subcategory of $\bold{TEmb}_{G}({X})$.
	\begin{itemize}
	\item[{\bf ($\pmb{\Pi}$)}]
	We say that $\mathcal{S}_{G}$ {\em contains all products} if given any pair of $G$-equivariant closed embeddings $\imath:{X}\mono Y_{\Sigma}$ and $\imath':{X}\mono Y_{\Sigma'}$  in $\mathcal{S}_{G}$, the composite $G$-equivariant embedding
	$$
	{X}\ \!\xymatrix{{}\ar@{^{(}->}^{\imath\times\imath'}[r]&{}}\!\ Y_{\Sigma}\times_{{K}}Y_{\Sigma'}\xrightarrow{\ \sim\ }Y_{\Sigma\times\Sigma'}
	$$
is an object in $\mathcal{S}_{G}$, and if each of the $G$-equivariant projections
	$
	Y_{\Sigma}\xleftarrow{\ \text{pr}_{1}\ }Y_{\Sigma\times\Sigma'}\xrightarrow{\ \text{pr}_{2}}Y_{\Sigma'}
	$
is a morphism in $\mathcal{S}_{G}$.
	\item[{\bf ($\pmb{\star}$)}]\vskip .2cm
	We say that $\mathcal{S}_{G}$ {\em satisfies condition ($\star$)} if there exists an affine open cover $\{U_{1},\dots,U_{r}\}$ of ${X}$ such that for each $1\le i\le r$ and each regular function $f$ on $U_{i}$, there exists a $G$-equivariant closed embedding $\imath:{X}\mono Y_{\Sigma}$ such that $U_{i}$ is the inverse image of a torus invariant open subset of $Y_{\Sigma}$ and $f$ is the pullback of a monomial on the dense torus in $Y_{\Sigma}$.
	\end{itemize}
\end{definition}

\begin{construction}
\normalfont
	The category of topological spaces equipped with a $G^{\op}$-action, i.e., the category of {\em $G^{\op}$-spaces}, admits inverse limits. If $\mathcal{S}_{G}$ is a system of $G$-equivariant toric embeddings of ${X}$, then we can form the inverse limit
	$$
	\varprojlim_{\mathcal{S}_{G}}\text{Trop}_{G}({X},\imath)
	$$
in the category of $G^{\op}$-spaces. It comes with a continuous $G$-equivariant map
	\begin{equation}\label{equation: canonical G-equivariant map}
	\pi:X^{\an}\ \longrightarrow\ \varprojlim_{\mathcal{S}_{G}}\text{Trop}_{G}({X},\imath).
	\end{equation}
\end{construction}

\begin{remark}\label{remark: A2}
In \cite{Wlod}, J. W\l odarczyk shows that if $K_{0}$ is algebraically closed and $X_{0}$ is normal, then $X_{0}$ admits a closed embedding into a toric $K_{0}$-variety if and only if $X_{0}$ is an {\em A${}_{2}$-variety}, i.e., $X_{0}$ satisfies
	\begin{itemize}
	\item[{\bf (A${}_{\bold{2}}$)}\!\!\!\!\!]\vskip .1cm \ \ \ every pair of points $x,y\in X_{0}$ lies in some affine open subset of $X_{0}$.\vskip .2cm
	\end{itemize}
W\l odarczyk's result can be understood as a generalization of the Chevalley-Kleiman criterion \cite[\S IV.2, Theorem 3]{Kleiman:66} \cite[\S I.9, Theorem 9.1]{Hartshorne:70}. In \cite[\S 4]{FGP}, P. Gross, S. Payne, and the present author use techniques developed by W\l odarczyk in \cite[\S 4]{Wlod} to show that (if one ignores $G$-action throughout), then for any system of closed toric embeddings satisfying the non-$G$-equivariant versions of Conditions ($\Pi$) and ($\star$), the map \eqref{equation: canonical G-equivariant map} is a homeomorphism \cite[Theorem 1.1]{FGP}. Using details of W\l odarczyk's result, this allows one to show that for any variety over an algebraically closed non-Archimedean field, admitting a closed embedding into at least one toric variety over an algebraically closed non-Archimedean field, there are abundant of examples of systems of toric embeddings for which the non-equivariant version of \eqref{equation: canonical G-equivariant map} is a homeomorphism \cite[Theorem 1.2]{FGP}.

	Our first two results in the present paper reproduce this theorem in the Galois-equivariant setting.
\end{remark}

\begin{Th}\label{theorem: main A}
\normalfont
	For any system $\mathcal{S}_{G}$ of $G$-equivariant toric embeddings satisfying Conditions ($\Pi$) and ($\star$) of Definition \ref{defintion: G-equivariant system}, the map \eqref{equation: canonical G-equivariant map} is a $G$-equivariant homeomorphism.
\end{Th}

\begin{Th}\label{theorem: main B}
\normalfont
	Suppose that $K_{0}$ is Henselian and that $X_{0}$ is a $K_{0}$-variety such that either $X$ is quasiprojective or $K_{0}$ is perfect and ${X}$ admits a closed embedding into at least one toric variety. Then ${X}$ admits systems of $G$-equivariant toric embeddings satisfying Conditions ($\Pi$) and ($\star$). The category $\bold{TEmb}_{G}({X})$ of all $G$-equivariant toric embeddings is one such system.
\end{Th}

\begin{remark}
	Theorems \ref{theorem: main A} and \ref{theorem: main B} together imply that the $G^{\op}$-space given by ${X}^{\an}$ {\em with its $G^{\op}$-action} can be reconstructed from Galois equivariant tropicalizations of ${X}$, giving an affirmative answer to Question \ref{questions: recovering Berkovich analytification tropically}.
	
	Because the canonical continuous surjection ${X}^{\an}\!\!\epi{X}$ used in formulating the universal property of the Berkovich analytification \cite[Theorem 3.4.1]{Berkovich:90} is $G$-equivariant, Theorem \ref{theorem: main B} implies that the underlying topological space of the ${K}$-variety ${X}$ {\em with its $G$-action} is encoded, in its entirety, in the totality of Galois equivariant tropicalizations of ${X}$. Succinctly: all Galois-action on a torically embeddable variety over a Henselian field is tropical.
	
	Our third main result gives an affirmative answer to Question \ref{questions: recovering orbits tropically}.
\end{remark}

\begin{Th}\label{theorem: main C}
\normalfont
	Suppose that $K_{0}$ is Henselian and that $X_{0}$ is a $K_{0}$-variety such that either $X$ is quasiprojective or $K_{0}$ is perfect and ${X}$ admits a closed embedding into at least one toric variety. Then for each point $x\in X^{\an}$ with finite $G^{\op}$-orbit, there exists a $G$-equivariant closed embedding $\imath:{X}\mono Y_{\Sigma}$ into a $G^{\op}$-twisted toric variety $Y_{\Sigma}$ such that the composite
	$$
	G^{\op}x\mono {X}^{\an}\xrightarrow{\ \text{trop}\ }\text{Trop}_{G}({X},\imath)
	$$
is a $G$-equivariant isomorphism of the $G^{\op}$-set $G^{\op}x$ onto its image in $\text{Trop}_{G}({X},\imath)$.
\end{Th}

\begin{remark}
	Theorem \ref{theorem: main C} says that a finite analytic Galois-orbit in an algebraic variety can always be recovered in a single Galois-equivariant tropicalization. Our last main result provides the beginning of an answer to the corresponding question for a different kind of finite geometric structure on a variety:
\end{remark}

\begin{Th}\label{theorem: main D}
\normalfont
	Suppose that $K_{0}$ is Henselian and that ${X}$ is projective. Let $F$ be any field of characteristic-$0$. Then for each $G$-equivariant closed embedding $\imath:X\mono Y_{\Sigma}$, we obtain an induced representation 
	$$
	G\ \ \lact\ \ H^{p,q}_{\text{trop}}\big(\ \!\text{Trop}({X},\imath),F\ \!\big)
	\ \ \ \ \mbox{and}\ \ \ \ 
	G^{\op}\ \ \lact\ \ H_{p,q}^{\text{trop}}\big(\ \!\text{Trop}({X},\imath),F\ \!\big)
	\ \ \ \mbox{for each}\ \ \ p,q\ge0
	$$
on the tropical cohomology and homology groups of Itenberg, Katzarkov, Mikhalkin, and Zharkov.
\end{Th}

\begin{remark}
	Because the proofs of Theorems \ref{theorem: main B} and \ref{theorem: main C} provide us with a large supply of $G$-equivariant closed embeddings of $X$, Theorem \ref{theorem: main D} provides us with a large supply of nontrivial, Artin Galois-representations in $F$-vector spaces.
\end{remark}
	
\end{subsection}


\begin{subsection}{Outline of the paper}
	In \S\ref{section: Galois equivariant tropicalizaiton} we review the theories of Galois-twisted and arithmetic toric varieties, we use them to define extended Galois-equivariant tropicalizations, and we prove Theorem  \ref{theorem: main D}. In \S\ref{section: Cox rings, characteristic spaces, and Galois-equivariant embeddings} we show how to use either a lemma of Payne or algorithms of W\l odarczyk or Berchtold and Hausen to construct a large supply of Galois-equivariant toric embeddings. In \S\ref{section: proofs of the main theorems} we prove Theorems \ref{theorem: main A}, \ref{theorem: main B}, and \ref{theorem: main C}.
\end{subsection}


\begin{subsection}{Acknowledgements}
	First, I thank the organizers of the Georgia Algebraic Geometry Symposium 2015, where the early ideas for this paper took shape. Thank you to Efrat Bank, Max Hully, and Dhruv Ranganathan for listening to and commenting on these ideas, and to Kristin Shaw and Frank Sottile for helpful conversations. Alexander Duncan played a crucial role by alerting me to his paper \cite{Duncan} and to the theory of twisted toric varieties. I thank him and Sam Payne for their enthusiasm and for detailed comments on drafts of the paper. Finally, I thank the organizers of the Tropical Geometry Seminar at l'Institut Mathématiques de Jussieu, namely Erwan Brugall\'e, Penka Georgieva, and Ilia Itenberg, for a speaking opportunity that greatly improved the paper.
	
	Research for the paper was conducted at University of Michigan, at l'Institut des Hautes \'Etudes Scientifiques, and at l'Institut Henri Poincar\'e. I thank all three institutions for their hospitality. Support for the author came from NSF RTG grant DMS-0943832 and from Laboratoire d'Excellence CARMIN.
\end{subsection}

\end{section}


\vskip .5cm
\begin{section}{Galois-equivariant tropicalization}\label{section: Galois equivariant tropicalizaiton}

	Having defined Galois-equivariant tropicalizations of closed subvarieties of Galois-twisted tori in \S\ref{subsection: Galois actions on tropicalizations} of the Introduction, we now extend this construction to provide Galois-equivariant tropicalizations of closed Galois-invariant subvarieties of all Galois-twisted toric varieties.


\begin{subsection}{Toric varieties over non-closed fields}
	We begin by providing a brief review the theory of toric varieties twisted by a Galois action, as developed by E. J. Elizondo, P. Lima-Filho, F. Sottile, and Z. Teitler in \cite{ELFST}, and by A. Duncan in \cite{Duncan}. This theory has a long gestation period that precedes \cite{ELFST} and \cite{Duncan}. A thorough account of this development appears in \cite[\S1]{ELFST}.

	Fix a field $K_{0}$, not necessarily equipped with a non-Archimedean absolute value. Fix a separable closure ${K}=K^{\text{sep}}_{0}$. If $Y_{0}$ is a $K_{0}$-variety, with $Y\Def Y_{0,{K}}$, then the $G$-action $G\ \lact\ {K}$ induces the {\em canonical $G^{\op}$-action} $G^{\op}\ \lact\ Y$, wherein the action of each $g\in G$ is given by the Cartesian diagram
	\begin{equation}\label{equation: Cartesian canonical}
	\begin{aligned}
	\xymatrix{
	Y
	\ar[r]^{g^{\ast}}
	\ar[d]
	&
	Y
	\ar[d]
	\\
	\text{Spec}_{\ \!}{{K}}
	\ar[r]_{g^{\ast}}
	&
	\text{Spec}_{\ \!}{{K}}.
	}
	\end{aligned}
	\end{equation}
A  {\em twisted $G^{\op}$-action on $Y_{}$} is any $G^{\op}$-action $G^{\op}\ \lact\ Y$ induced by the action $H^{\op}\ \lact\ Y_{0,L}$ of some finite quotient $H=G/Z\cong\text{Gal}(L/K_{0})$ with fixed field $L={K}^{Z}$, such that for each $h\in H$ the diagram
	$$
	\xymatrix{
	Y_{0,L}
	\ar[r]^{h^{\ast}}
	\ar[d]
	&
	Y_{0,L}
	\ar[d]
	\\
	\text{Spec}_{\ \!}{L}
	\ar[r]_{h^{\ast}}
	&
	\text{Spec}_{\ \!}{L}
	}
	$$
commutes. For any twisted $G^{\op}$-action on $Y$, the diagram \eqref{equation: Cartesian canonical} commutes for each $g\in G$, although the diagram is not Cartesian in general. A {\em $K_{0}$-form} of $Y$ is any $K_{0}$-variety $Y'_{0}$ that admits an isomorphism of ${K}$-varieties
	$$
	Y'_{0,{K}}
	\ \cong\ 
	Y.
	$$

\begin{definition}[{\cite[\S3]{ELFST}}]\label{definition: arithmetic toric variety}
\normalfont
	An {\em arithmetic toric variety over $K_{0}$} is a triple $(Y_{0},T_{0},\mu)$ consisting of:
	\begin{itemize}
	\item[{\bf (i)}]\vskip .1cm
	a $K_{0}$-variety $Y_{0}$;
	\item[{\bf (ii)}]\vskip .2cm
	a (not-necessarily split) $K_{0}$-torus $T_{0}$;\footnote{See \S\ref{subsection: Galois actions on tropicalizations}}
	\item[{\bf (iii)}]\vskip .2cm
	a faithful $T_{0}$-action $\mu:T_{0}\ \lact\ Y_{0}$ with dense orbit in $Y_{0}$. \vskip .2cm
	\end{itemize}
We often write an arithmetic toric variety as simply $Y_{0}$, leaving $T_{0}$ and $\mu$ implicit.
\end{definition}

\begin{remark}
	Fix a split $K_{0}$-torus $T_{0}=\text{Spec}_{\ \!}K_{0}[M]$ with character lattice $M$. Define $T\Def T_{0,{K}}$. Equip $M$ with the trivial $G$-action. Then the automorphism group of $T$ {\em as a ${K}$-scheme} is
	$$
	\text{Aut}_{{K}}(T)
	\ \cong\ 
	T({K})\rtimes\text{GL}_{\ZZ}(M)^{\op}.
	$$
It comes with a $G^{\op}$-action that restricts to the canonical action on the first factor $T({K})$ and to the trivial action on the second factor $\text{GL}_{\ZZ}(M)^{\op}$. 

	Let $T_{0}\ \lact\ Y_{0}$ be an arithmetic toric variety over $K_{0}$. Then ${Y}=Y_{0,{K}}$ is isomorphic to a toric ${K}$-variety $Y_{\Sigma}$ for some fan $\Sigma$ in $N_{\RR}$. For each $g\in G$, the resulting morphism {\em of $K_{0}$-schemes}
	$$
	g^{\ast}:Y_{\Sigma}\longrightarrow Y_{\Sigma}
	$$
is $T$-equivariant. Up to translation by $T({K})$, it is a toric morphism, which is to say that up to translation by some element $y_{g}\in T({K})$, it comes from a $\ZZ$-integral automorphism $\varphi_{g}\ \lact\ N_{\RR}\Def\text{Hom}_{\ZZ}(M,\RR)$ that leaves the fan $\Sigma$ invariant. In this way, the twisted $G^{\op}$-action $G^{\op}\ \lact\ Y_{\Sigma}$ gives rise to a map of sets
	\begin{equation}\label{equation: induced cocycle}
	(y_{\sm{(-)}},\varphi_{\sm{(-)}})\ \!:\ G^{\op}\longrightarrow T({K})\rtimes\text{Aut}(\Sigma),
	\end{equation}
where
	$$
	\text{Aut}(\Sigma)
	\ \Def\ 
	\big\{\psi\in\text{GL}_{\ZZ}(N)\ \!:\ \Sigma\mbox{ is invariant under }N_{\RR}\xrightarrow{\ \psi\ }N_{\RR}\big\}.
	$$
Equip $\text{Aut}(\Sigma)$ with the trivial $G$-action, ${T}({K})$ with its canonical $G$-action, and ${T}({K})\rtimes\text{Aut}(\Sigma)$ with the resulting product $G$-action.
\end{remark}

\begin{proposition}[{\cite[\S III.1, Proposition 5]{Serre:2002}}, {\cite[Lemma 3.1 \& Theorem 3.2]{ELFST}}]\label{proposition: classification of geometrically projective arithmetic toric varieties}
\normalfont
	{\bf (i).} The set of $G^{\op}$-twisted actions on a ${K}$-variety $Y$, up to $G$-equivariant automorphisms, is in natural bijection with the non-Abelian Galois cohomology set
	\begin{equation}\label{equation: classifying cohomology group}
	H^{1}\big(G^{\op},\text{Aut}_{{K}}(Y)\big),
	\end{equation}
where $\text{Aut}_{{K}}(Y)$ denotes the group of automorphisms of $Y$ over ${K}$. If $Y$ is quasiprojective, then the set \eqref{equation: classifying cohomology group} is also in bijection with the set of twisted $K_{0}$-forms of $Y$.

	{\bf (ii).} The map \eqref{equation: induced cocycle} is a $1$-cocycle on $G^{\op}$. If $Y$ is a quasiprojective toric ${K}$-variety with fan $\Sigma$, then the construction of the $1$-cocycle \eqref{equation: induced cocycle} gives rise to a bijection
	$$
	\begin{aligned}
	\left.
	\left\{\!\!
	\begin{array}{c}
	\mbox{arithmetic toric $K_{0}$-varieties $T'_{0}\ \lact\ Y_{0}$}
	\\
	\mbox{admitting ${K}$-variety isomorphisms}
	\\
	\mbox{$T'\cong T$ and $Y\cong Y_{\Sigma}$ s.t. the diagram}
	\\
	\begin{aligned}
	\xymatrix{
	T'\times_{K}Y
	\ar[r]
	\ar[d]_{\wr}
	&
	Y
	\ar[d]^{\wr}
	\\
	T\times_{K} Y_{\Sigma}
	\ar[r]
	&
	Y_{\Sigma}
	}
	\end{aligned}
	\mbox{\ \ commutes}
	\end{array}
	\!\!\right\}
	\right/
	\mbox{iso.}
	\ \ \ \xrightarrow{\ \ \ \mbox{$\sim$}\ \ \ }\ \ \ 
	H^{1}\big(G^{\op},{T}({K})\!\rtimes\!\text{Aut}(\Sigma)\big).
	\end{aligned}
	$$
\end{proposition}

\begin{corollary}\label{corollary: arithmetic toric varieties from homomorphisms}
\normalfont
	{\bf (i).} Given any lattice $N$ and any projective fan $\Sigma$ in $N_{\RR}$, each group homomorphism
	\begin{equation}\label{equation: morphism cocycle}
	\varphi:G^{\op}\longrightarrow\text{Aut}(\Sigma)
	\end{equation}
determines a $1$-cocyle \eqref{equation: induced cocycle} with $y_{g}$ the identity in $T({K})$ for all $g\in G$. It therefore defines a toric ${K}$-variety $Y_{\Sigma}$ equipped with a twisted $G^{\op}$-action $G^{\op}\ \lact\ Y_{\Sigma}$. The representation of $G$ in finite rank $\ZZ$-lattices
	$$
	G\xrightarrow{\ \varphi^{\text{op}}\ }\text{Aut}(\Sigma)^{\op}\mono\text{Aut}(M)
	$$
is exactly the image of the corresponding $K_{0}$-twisted form of $T$ under the equivalence \eqref{equation: functor on free lattices}.
	
	{\bf (ii).} When $Y$ is quasiprojective, each homomorphism \eqref{equation: morphism cocycle}
has a uniquely associated arithmetic toric variety $\mu_{\varphi}:T_{0}\ \lact\ Y_{0}$ over $K_{0}$.
\hfill
$\square$
\end{corollary}

\begin{definition}\label{definition: twisted toric variety}
\normalfont
	A {\em $G^{\op}$-twisted toric variety over ${K}$} is any toric ${K}$-variety $Y_{\Sigma}$ equipped with a twisted $G^{\op}$-action arising from a homomorphism \eqref{equation: morphism cocycle} as in Corollary \ref{corollary: arithmetic toric varieties from homomorphisms}.(i). If $G^{\op}\ \lact\  Y_{\Sigma}$ and $G^{\op}\ \lact\ Y_{\Sigma'}$ are $G^{\op}$-twisted toric varieties, then a {\em morphism} $Y_{\Sigma}\longrightarrow Y_{\Sigma'}$ of $G^{\op}$-twisted toric varieties is any $G$-equivariant morphism of toric ${K}$-varieties induced by a $G$-equivariant morphism of fans $\Sigma\longrightarrow\Sigma'$.
\end{definition}

\begin{remark}
	We can construct the $G^{\op}$-twisted toric variety $G^{\op}\ \lact\ Y_{\Sigma}$ in Corollary \ref{corollary: arithmetic toric varieties from homomorphisms} explicitly. Observe that for each element $g\in G$ and each cone $\sigma\in \Sigma$, the fact that the actions $g\ \lact\ M$ and $g^{\ast}\ \lact\ N$ are adjoint implies that
	$$
	\begin{array}{rcl}
	g(S_{g^{\ast}(\sigma)})
	&
	\Def
	&
	\big\{
	\ \!g(u)\in M\ \!:\ \!\big\langle u,g^{\ast}(v)\big\rangle\ge 0\ \mbox{for all}\ g^{\ast}(v)\in g^{\ast}(\sigma)
	\big\}
	\\[8pt]
	&
	=
	&
	\big\{
	\ \!g(u)\in M\ \!:\ \!\big\langle g(u),v\big\rangle\ge 0\ \mbox{for all}\ v\in \sigma
	\big\}
	\\[8pt]
	&
	=
	&
	S_{\sigma}.
	\end{array}
	$$
Thus for each $g\in G$, the automorphism $g:M\xrightarrow{\ \sim\ }M$ restricts to an isomorphism of semigroups
	\begin{equation}\label{equation: semigroup homomorphism}
	g:S_{g^{\ast}(\sigma)}\xrightarrow{\ \sim\ }S_{\sigma}.
	\end{equation}
Letting $g$ act on ${K}$, this induces a $K_{0}$-algebra (!) isomorphism
	\begin{equation}\label{equation: action of g on semigroup algebras}
	g:{K}[S_{g^{\ast}(\sigma)}]\xrightarrow{\ \sim\ }{K}[S_{\sigma}],
	\end{equation}
hence a commutative diagram of $K_{0}$-scheme morphisms
	\begin{equation}\label{equation: g switching cone morphism}
	\begin{aligned}
	\xymatrix{
	Y_{\sigma}
	\ar[r]^{g^{\ast}\ \ \ }_{\sim\ \ \ \ \ }
	\ar@{->>}[d]
	&
	Y_{g^{\ast}(\sigma)}
	\ar@{->>}[d]
	\\
	\text{Spec}_{\ \!}{K}
	\ar[r]_{g^{\ast}}^{\sim}
	&
	\text{Spec}_{\ \!}{K}.\!
	}
	\end{aligned}
	\end{equation}
It is easy to see that for a given $g\in G$, the diagrams \eqref{equation: g switching cone morphism} commute with the inclusions of toric affine opens induced by inclusions of cones $\tau\subset \sigma$, and thus glue to produce $Y_{\Sigma}$ with twisted $G^{\op}$-action.
\end{remark}

\begin{example}\label{example: Brauer-Severi A}
\normalfont
{\bf Brauer-Severi variety arising as an arithmetic toric variety.}
	Let $K_{0}=\mathbb{C}(\!(t)\!)$, so that $K=\mathbb{C}(\!(t)\!)^{\text{sep}}\cong\CC(\!(t^{\QQ})\!)$, the field of Puiseux series, and $G=\text{Gal}\big(\CC(\!(t^{\QQ})\!)\big/\CC(\!(t)\!)\big)\cong\widehat{\ZZ}$. Let $Y=\PP^{2}_{\sm{\mathbb{C}(\!(t)\!)}}$. Define $N\Def\ZZ^{2}$, and let $\Sigma\subset$ be the standard fan in $\RR^{2}$ that realizes $\PP^{2}_{\sm{\CC(\!(t^{\QQ})\!)}}$ as a toric $\CC(\!(t^{\QQ})\!)$-variety, with rays generated by the vectors $e_{1}=(1,0)$, $e_{2}=(0,1)$, and $(-1,-1)$. Consider the 3-periodic matrix
	\begin{equation}\label{equation: 3-periodic matrix}
	\left(\!
	\begin{array}{cc}
	0 & \!-1
	\\[2pt]
	1 & \!-1
	\end{array}
	\!\right)
	\ \in\ 
	\text{SL}_{2}(\ZZ)
	\end{equation}
acting on $N_{\RR}$. As pictured in Figure \ref{figure: example of action}, the fan $\Sigma$ is invariant with respect to the action of this matrix. The matrix provides a cyclic permutation of the cones and rays in $\Sigma$.
	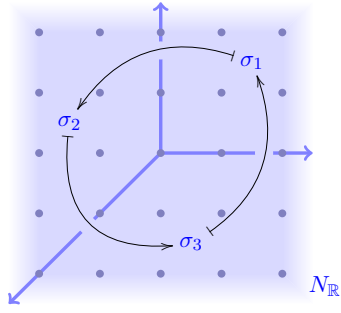
\begin{figure}[h!]
	$$
	\scalebox{.8}{$
	\begin{xy}
	(0,0)*+{
	\begin{tikzpicture}
	\fill[blue!15, path fading=north] (-2,2) -- (2,2) -- (2.5,2.5) -- (-2.5,2.5);
	\fill[blue!15, path fading=south] (-2,-2) -- (2,-2) -- (2.5,-2.5) -- (-2.5,-2.5);
	\fill[blue!15, path fading=east] (2,2) -- (2,-2) -- (2.5,-2.5) -- (2.5,2.5);
	\fill[blue!15, path fading=west] (-2,2) -- (-2,-2) -- (-2.5,-2.5) -- (-2.5,2.5);
	\fill[blue!15] (-2,-2) -- (2,-2) -- (2,2) -- (-2,2);
	\draw[blue!50, ultra thick, ->] (0,0) -- (2.5,0);
	\draw[blue!50, ultra thick, ->] (0,0) -- (0,2.5);
	\draw[blue!50, ultra thick, ->] (0,0) -- (-2.5,-2.5);
	\fill[blue!15] (1.7,0) circle (.15);
	\fill[blue!15] (0,1.7) circle (.15);
	\fill[blue!15] (-1.2,-1.2) circle (.15);
	\fill[black!50!blue!50] (0,0) circle (.065);
	\fill[black!50!blue!50] (1,0) circle (.065);
	\fill[black!50!blue!50] (2,0) circle (.065);
	\fill[black!50!blue!50] (-1,0) circle (.065);
	\fill[black!50!blue!50] (-2,0) circle (.065);
	\fill[black!50!blue!50] (0,1) circle (.065);
	\fill[black!50!blue!50] (1,1) circle (.065);
	\fill[black!50!blue!50] (2,1) circle (.065);
	\fill[black!50!blue!50] (-1,1) circle (.065);
	\fill[black!50!blue!50] (-2,1) circle (.065);
	\fill[black!50!blue!50] (0,2) circle (.065);
	\fill[black!50!blue!50] (1,2) circle (.065);
	\fill[black!50!blue!50] (2,2) circle (.065);
	\fill[black!50!blue!50] (-1,2) circle (.065);
	\fill[black!50!blue!50] (-2,2) circle (.065);
	\fill[black!50!blue!50] (0,-1) circle (.065);
	\fill[black!50!blue!50] (1,-1) circle (.065);
	\fill[black!50!blue!50] (2,-1) circle (.065);
	\fill[black!50!blue!50] (-1,-1) circle (.065);
	\fill[black!50!blue!50] (-2,-1) circle (.065);
	\fill[black!50!blue!50] (0,-2) circle (.065);
	\fill[black!50!blue!50] (1,-2) circle (.065);
	\fill[black!50!blue!50] (2,-2) circle (.065);
	\fill[black!50!blue!50] (-1,-2) circle (.065);
	\fill[black!50!blue!50] (-2,-2) circle (.065);
	\end{tikzpicture}
	};
	(15,15)*+{{\color{blue} \sigma_{1}}}="1";
	(-15,5)*+{{\color{blue} \sigma_{2}}}="2";
	(5,-15)*+{{\color{blue} \sigma_{3}}}="3";
	(27,-22)*+{{\color{blue} N_{\RR}}};
	{\ar@{|->}@/_20pt/ "1"; "2"};
	{\ar@{|->}@/_27pt/ "2"; "3"};
	{\ar@{|->}@/_20pt/ "3"; "1"};
	\end{xy}
	$}
	$$
\caption{
The action of the $3$-periodic matrix \eqref{equation: 3-periodic matrix} on the fan $\Sigma$.
}
\label{figure: example of action}
\end{figure}

\noindent
Fix an isomorphism
	$$
	\ZZ/3\ZZ
	\ \cong\ 
	\text{Gal}\big(\CC(\!(t^{\frac{1}{3}})\!)\big/\CC(\!(t)\!)\big)^{\op}.
	$$
Let $\ZZ/3\ZZ\longrightarrow\text{Aut}(\Sigma)$ be the group homomorphism taking $1\mapsto\left(\begin{smallmatrix} 0 & -1 \\ 1 & -1 \end{smallmatrix}\right)$, and consider the composite
	$$
	\varphi\ :\ \ G^{\op}\epi \ZZ/3\ZZ\longrightarrow\text{Aut}(\Sigma).
	$$
By Corollary \ref{corollary: arithmetic toric varieties from homomorphisms}.(ii), the homomorphism $\varphi$ determines an arithmetic toric variety $\mu:T_{0}\ \lact\ Y_{0}$ over $\CC(\!(t)\!)$. The $\CC(\!(t)\!)$-variety $Y_{0}$ is a $2$-dimensional Brauer-Severi variety \cite{Artin:82} that comes equipped with the faithful, dense action of a non-split $2$-dimensional torus $T_{0}$.
\end{example}

\end{subsection}


\begin{subsection}{Galois equivariant extended tropicalizations of twisted toric varieties}
	Equip $K_{0}$ with the structure of a (rank-$1$) non-Archimedean field, i.e., equip $K_{0}$ with a non-Archimedean absolute value
	$$
	|-|:K_{0}\longrightarrow\RR_{\ge0}.
	$$
Fix a separable closure $K=K^{\text{sep}}_{0}$, and let $|-|:{K}\longrightarrow\RR_{\ge0}$ be a non-Archimedean valuation extending our chosen valuation on $K_{0}$. When $K_{0}$ is Henselian, for instance when $K_{0}$ is complete with respect to its valuation, the extension of $|-|$ to ${K}$ is unique \cite[Lemma 4.1.1]{EP:05}. Define
	$$
	G
	\ \Def\ 
	\text{Gal}(K/K_{0}).
	$$
	
	Fix a $G^{\op}$-twisted toric variety $G^{\op}\ \lact\ Y_{\Sigma}$ determined by a group homomorphism $\varphi:G^{\op}\longrightarrow\text{Aut}(\Sigma)$ for some fan $\Sigma\subset N_{\RR}$ as in Corollary \ref{corollary: arithmetic toric varieties from homomorphisms}.(i). Let us temporarily ignore the Galois-structure on $Y_{\Sigma}$, interpreting $Y_{\Sigma}$ as a toric ${K}$-variety without $G^{\op}$-action. S. Payne \cite[\S3]{P09} explains how to associate a Hausdorff space to $Y_{\Sigma}$, called the extended tropicalization of $Y_{\Sigma}$. In detail, each cone $\sigma\in\Sigma$ has a corresponding topological space
	$$
	\text{Trop}(Y_{\sigma})
	\ \Def\ 
	\text{Hom}_{\bold{SGrp}}\big(S_{\sigma},\RR\sqcup\{\infty\}\big),
	$$
where $S_{\sigma}$ is the semigroup $S_{\sigma}\Def\sigma^{\vee}\cap M$ in $M_{\RR}$, and where $\RR\sqcup\{\infty\}$ is the ordered semigroup with $\infty>r$ for all $r\in\RR$, equipped with its order topology. Inclusions of cones $\tau\subset\sigma$ in $\Sigma$ induce open embeddings of topological spaces
	$$
	\text{Trop}(Y_{\tau})\ \mono\text{Trop}(Y_{\sigma}).
	$$
The {\em extended tropicalization of $Y_{\Sigma}$}, denoted $\text{Trop}(Y_{\Sigma})$, is the topological space obtained by gluing along these inclusions:
	$$
	\text{Trop}(Y_{\Sigma})
	\ \Def\ 
	\varinjlim_{\sigma\in\Sigma}\text{Trop}(Y_{\sigma}).
	$$
It contains $N_{\RR}=\text{Trop}(Y_{\{0\}})$ as a dense open subset. The extended tropicalization comes with a continuous surjective {\em tropicalization map}
	\begin{equation}\label{equation: extended tropicalization map}
	\text{trop}:Y^{\an}_{\Sigma}\epi\text{Trop}(Y_{\Sigma}).
	\end{equation}
On each $Y_{\sigma}=\text{Spec}_{\ \!}{K}[S_{\sigma}]$, the map \eqref{equation: extended tropicalization map} takes a point $y\in Y^{\an}_{\sigma}$, given by a multiplicative seminorm $|-|_{y}:{K}[S_{\sigma}]\longrightarrow\RR_{\ge0}$, and sends it to the semigroup homomorphism
	$$
	S_{\sigma}\mono{K}[S_{\sigma}]\xrightarrow{\ |-|_{y}\ }\RR_{\ge0}\xrightarrow{-\log\ \ }\RR\sqcup\{\infty\}.
	$$
	
\begin{lemma}\label{lemma: extended G-action}
\normalfont
	The $G^{\op}$-action $G^{\op}\ \lact\ N_{\RR}$ extends to a $G^{\op}$-action
	\begin{equation}\label{equation: induced Galois-action on tropicalization - first appearance}
	G^{\op}\ \lact\ \text{Trop}(Y_{\Sigma}).
	\end{equation}
If $K_{0}$ is Henselian, then the twisted $G^{\op}$-action $G^{\op}\ \lact\ Y_{\Sigma}$ induces a $G^{op}$-action
	\begin{equation}\label{equation: Berkovich Galois action}
	G^{\op}\ \lact\ Y^{\an}_{\Sigma},
	\end{equation}
and the tropicalization map \eqref{equation: extended tropicalization map} is $G$-equivariant with respect to the $G^{\op}$-actions \eqref{equation: induced Galois-action on tropicalization - first appearance} and \eqref{equation: Berkovich Galois action}.
\end{lemma}
\begin{proof}
	For $g\in G$, the semigroup homomorphism \eqref{equation: semigroup homomorphism} induces a homeomorphism of topological spaces
	\begin{equation}\label{equation: local homeomorphism from cone}
	g^{\ast}:\text{Trop}(Y_{\sigma})\xrightarrow{\ \sim\ }\text{Trop}(Y_{g^{\ast\!}\sigma})
	\end{equation}
for each $\sigma\in\Sigma$. Because the semigroup homomorphisms \eqref{equation: semigroup homomorphism} commute with the morphisms dual to inclusions of cones $\tau\subset\sigma$, the homeomorphisms \eqref{equation: local homeomorphism from cone} glue to give an automorphism
	$$
	g^{\ast}\ \lact\ \text{Trop}(Y_{\Sigma}).
	$$
Taking these maps for all $g\in G$, we obtain a $G^{\op}$-action on $\text{Trop}(Y_{\Sigma})$.

	For each $\sigma\in\Sigma$, a point $y\in Y^{\text{an}}_{\sigma}$ is given by a multiplicative seminorm $|-|_{y}:{K}[S_{\sigma}]\longrightarrow\RR_{\ge0}$ who's restriction to ${K}\subset {K}[S_{\sigma}]$ is our fixed absolute value $|-|$ on ${K}$. If $K_{0}$ is Henselian, then for each $g\in G$, the $K_{0}$-algebra isomorphism \eqref{equation: action of g on semigroup algebras} fits into a commutative diagram
	\begin{equation}\label{equation: commuting action on semigroup algebras}
	\begin{aligned}
	\begin{xy}
	(-20,0)*+{\ {K}\ \!}="-1";
	(-20,-16)*+{\ \!\ {K}\ }="0";
	(0,0)*+{{K}[S_{g^{\ast\!}\sigma}]}="1";
	(0,-16)*+{{K}[S_{\sigma}]}="2";
	(25,-8)*+{\RR_{\ge0}.}="3";
	(15,5)*+{\mbox{{\smaller\smaller\smaller\smaller\smaller $|\!-\!|$}}};
	(15,-21)*+{\mbox{{\smaller\smaller\smaller\smaller\smaller $|\!-\!|$}}};
	{\ar_{g} "-1"; "0"};
	{\ar@{_{(}->} "-1"; "1"};
	{\ar@{^{(}->} "0"; "2"};
	{\ar_{g} "1"; "2"};
	{\ar@/_5pt/^{|g(-)|_{y}\ \ } "1"; "3"};
	{\ar@/^5pt/_{|-|_{y}\!\!\!\!\!\!} "2"; "3"};
	{\ar@/^30pt/ "-1"; "3"};
	{\ar@/_30pt/ "0"; "3"};
	\end{xy}
	\end{aligned}
	\end{equation}
In particular, the composite multiplicative seminorm $|-|_{g^{\ast}(y)}\Def\big|g(-)\big|_{v}$ defines a point $g^{\ast}(y)\in Y^{\an}_{g^{\ast}(\sigma)}$. The resulting  homeomorphisms $g^{\ast}:Y^{\an}_{\sigma}\xrightarrow{\ \sim\ }Y^{\an}_{g^{\ast}(\sigma)}$ glue to give an automorphism $g^{\ast}\ \lact\ Y^{\an}_{\Sigma}$, and together these give a $G^{\op}$-action \eqref{equation: Berkovich Galois action}. Finally, commutativity of the diagrams
	\begin{equation}\label{equation: commuting needle point}
	\begin{aligned}
	\begin{xy}
	(-20,0)*+{S_{g^{\ast\!}(\sigma)}}="-1";
	(-20,-16)*+{S_{\sigma}}="0";
	(0,0)*+{{K}[S_{g^{\ast\!}(\sigma)}]}="1";
	(0,-16)*+{{K}[S_{\sigma}]}="2";
	(25,-8)*+{\RR_{\ge0}}="3";
	(47.5,-8)*+{\RR\sqcup\{\infty\}}="4";
	{\ar_{g} "-1"; "0"};
	{\ar@{^{(}->} "-1"; "1"};
	{\ar@{^{(}->} "0"; "2"};
	{\ar_{g} "1"; "2"};
	{\ar@/_5pt/^{|-|_{g^{\ast\!}(y)}\!\!\!\!\!\!} "1"; "3"};
	{\ar@/^5pt/_{|-|_{y}\ \ } "2"; "3"};
	{\ar^{-\log\ \ \ \ \ \ } "3"; "4"};
	\end{xy}
	\end{aligned}
	\end{equation}
implies that the tropicalization map \eqref{equation: extended tropicalization map} is equivariant with respect to the actions constructed.
\end{proof}
	
\begin{definition}
\normalfont
	If $Y_{\Sigma}$ is a $G^{\op}$-twisted toric ${K}$-variety, then the {\em $G$-equivariant extended tropicalization} of $Y_{\Sigma}$ is the extended tropicalization $\text{Trop}(Y_{\Sigma})$ equipped with the $G^{\op}$-action constructed in the proof of Lemma \ref{lemma: extended G-action}. We denote it 
	$$
	\text{Trop}_{G}(Y_{\Sigma}).
	$$
When $K_{0}$ is Henselian, we refer to \eqref{equation: extended tropicalization map} as the {\em $G^{\op}$-equivariant tropicalizaiton map}.
\end{definition}

\begin{example}\label{example: Brauer-Severi B}
\normalfont
	Let $K_{0}=\CC(\!(t)\!)$, with $Y_{0}$ the $2$-dimensional Brauer-Severi $\CC(\!(t)\!)$-variety realized as an arithmetic toric variety $\mu:T_{0}\ \lact\ Y_{0}$ in Example \ref{example: Brauer-Severi A}. The compactifications $\overline{\sigma_{i}}$ of each of the three $2$-dimensional cones $\sigma_{i}$ in $\Sigma$ fit together to form $\text{Trop}(Y_{\Sigma})$, which appears in Figure \ref{figure: example of tropical action}. The $\ZZ/3\ZZ$-action $\ZZ/3\ZZ\ \lact\ N_{\RR}$, cyclically permuting the cones $\sigma_{i}$, $1\le i\le 3$, extends to a $\ZZ/3\ZZ$-action on the whole space $\text{Trop}(Y_{\Sigma})$. It cyclically permutes the three compactified cones $\overline{\sigma_{i}}$ as depicted in Figure \ref{figure: example of tropical action}. Precomposing with $G\epi \ZZ/3\ZZ$, we get a $G^{\op}$-action on $\text{Trop}(Y_{\Sigma})$.
	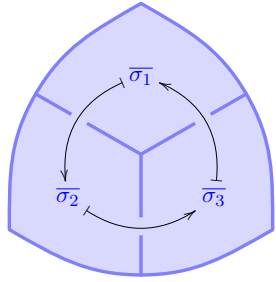
\begin{figure}[h!]
	$$
	\scalebox{.8}{$
	\begin{xy}
	(0,0)*+{
	\begin{tikzpicture}
	\fill[blue!15, ultra thick, line cap=round, rounded corners=.1mm] (0,0) -- ({2*cos(30)},{2*sin(30)}) to [out=120,in=-30] ({{2.5*cos(90)}},{{2.5*sin(90)}}) to [out=-150,in=60] ({{2*cos(150)}},{{2*sin(150)}}) -- (0,0);
	\fill[blue!15, ultra thick, line cap=round, rounded corners=.1mm] (0,0) -- ({2*cos(30+120)},{2*sin(30+120)}) to [out=-120,in=90] ({{2.5*cos(90+120)}},{{2.5*sin(90+120)}}) to [out=-30,in=180] ({{2*cos(150+120)}},{{2*sin(150+120)}}) -- (0,0);
	\fill[blue!15, ultra thick, line cap=round, rounded corners=.1mm] (0,0) -- ({2*cos(30-120)},{2*sin(30-120)}) to [out=0,in=-150] ({{2.5*cos(90-120)}},{{2.5*sin(90-120)}}) to [out=90,in=-60] ({{2*cos(150-120)}},{{2*sin(150-120)}}) -- (0,0);
	\draw[blue!50, ultra thick, line cap=round, rounded corners=.1mm] ({2*cos(30)},{2*sin(30)}) to [out=120,in=-30] ({{2.5*cos(90)}},{{2.5*sin(90)}}) to [out=-150,in=60] ({{2*cos(150)}},{{2*sin(150)}});
	\draw[blue!50, ultra thick, line cap=round, rounded corners=.1mm] ({2*cos(30+120)},{2*sin(30+120)}) to [out=-120,in=90] ({{2.5*cos(90+120)}},{{2.5*sin(90+120)}}) to [out=-30,in=180] ({{2*cos(150+120)}},{{2*sin(150+120)}});
	\draw[blue!50, ultra thick, line cap=round, rounded corners=.1mm] ({2*cos(30-120)},{2*sin(30-120)}) to [out=0,in=-150] ({{2.5*cos(90-120)}},{{2.5*sin(90-120)}}) to [out=90,in=-60] ({{2*cos(150-120)}},{{2*sin(150-120)}});
	\draw[blue!50, ultra thick, line cap=round, rounded corners=.1mm] (0,0) -- ({2*cos(30)},{2*sin(30)});
	\draw[blue!50, ultra thick, line cap=round, rounded corners=.1mm] (0,0) -- ({2*cos(150)},{2*sin(150)});
	\draw[blue!50, ultra thick, line cap=round, rounded corners=.1mm] (0,0) -- ({2*cos(-90)},{2*sin(-90)});
	\fill[blue!15] ({1.175*cos(30)},{1.175*sin(30)}) circle (.15);
	\fill[blue!15] ({1.175*cos(30+120)},{1.175*sin(30+120)}) circle (.15);
	\fill[blue!15] ({1.195*cos(-90)},{1.195*sin(-90)}) circle (.15);
	\end{tikzpicture}
	};
	(0,10.5)*+{{\color{blue} \overline{\sigma_{1}}}}="1";
	(-12,-9.5)*+{{\color{blue} \overline{\sigma_{2}}}}="2";
	(12,-9.5)*+{{\color{blue} \overline{\sigma_{3}}}}="3";
	{\ar@{|->}@/_15pt/ "1"; "2"};
	{\ar@{|->}@/_15pt/ "2"; "3"};
	{\ar@{|->}@/_15pt/ "3"; "1"};
	\end{xy}
	$}
	$$
\caption{Action of $\ZZ/3\ZZ$ on $\text{Trop}(Y_{\Sigma})$ from Example \ref{example: Brauer-Severi A}.}
\label{figure: example of tropical action}
\end{figure}
\end{example}

\begin{ssubsection}
\normalfont
{\bf Extended tropicalizations of $\pmb{K}$-varieties with canonical $\pmb{G^{\op}}$-action.}
	Let $G^{\op}\ \lact\ Y_{\Sigma}$ be a $G^{\op}$-twisted toric variety over ${K}$. Let $X_{0}$ be a $K_{0}$-variety. Equip $X\Def X_{0,{K}}$ with its canonical $G^{\op}$-action. Fix an affine open cover $\{\text{Spec}_{\ \!}A_{i}\hookrightarrow X_{0}\}_{i\in I}$. If $K_{0}$ is Henselian, then for each $g\in G$ and each point $x\in(\text{Spec}_{\ \!}A_{i,{K}})^{\an}$, we have have a commutative diagram 
	$$
	\begin{xy}
	(-16,0)*+{\ {K}\ \!}="-1";
	(-16,-12)*+{\ \!\ {K}\ }="0";
	(0,0)*+{A_{i,{K}}\!\!\!}="1";
	(0,-12)*+{A_{i,{K}}\!\!\!}="2";
	(25,-6)*+{\RR_{\ge0},}="3";
	(20,2)*+{\mbox{{\smaller\smaller\smaller\smaller\smaller $|\!-\!|$}}};
	(20,-14)*+{\mbox{{\smaller\smaller\smaller\smaller\smaller $|\!-\!|$}}};
	{\ar_{g} "-1"; "0"};
	{\ar@{_{(}->} "-1"; "1"};
	{\ar@{^{(}->} "0"; "2"};
	{\ar_{g} "1"; "2"};
	{\ar@/_4pt/^{\mbox{{\smaller\smaller\smaller\smaller\smaller $|\!-\!|_{g^{\ast\!}(x)}$}}\ \ \ } "1"; "3"};
	{\ar@/^4pt/_{\mbox{{\smaller\smaller\smaller\smaller\smaller $|\!-\!|_{x}$}}\ } "2"; "3"};
	{\ar@/^21pt/ "-1"; "3"};
	{\ar@/_21pt/ "0"; "3"};
	\end{xy}
	$$
and thus an induced {\em canonical $G^{\op}$-action} $G^{\op}\ \lact\ X^{\an}$. Each $G$-equivariant closed embedding
	$$
	\imath:{X}\mono Y_{\Sigma}
	$$
induces a $G$-equivariant closed embedding of topological spaces $\imath^{\an}:{X}^{\an}\mono Y^{\an}_{\Sigma}$. The ({\em extended}) {\em tropicalization of ${X}$} ({\em with respect to $\imath$}), denoted $\text{Trop}({X},\imath)$, is the image
	$$
	\text{Trop}({X},\imath)
	\ \Def\ 
	\text{im}\big(\ {X}^{\an}\!\!\xrightarrow{\ \ \imath^{\an}\ \ }Y^{\an}_{\Sigma}\xrightarrow{\text{trop}}\text{Trop}(Y_{\Sigma})\ \big).
	$$
The restriction of \eqref{equation: extended tropicalization map} to $X^{\an}$ is a continuous surjective map
	\begin{equation}\label{equation: tropicalization map on X}
	\text{trop}:{X}^{\an}\epi\text{Trop}(X,\imath)
	\end{equation}
that we call the {\em tropicalization map on $X^{\text{an}}$} ({\em with respect to $\imath$}). 
\end{ssubsection}
	
\begin{lemma}\label{lemma: invariance and equivariance}
\normalfont
	If $K_{0}$ is Henselian, then the tropicalization $\text{Trop}({X},\imath)$ is a $G^{\op}$-invariant subspace of $\text{Trop}(Y_{\Sigma})$ under the $G^{\op}$-action in Lemma \ref{lemma: extended G-action}, and the tropicalization map \eqref{equation: tropicalization map on X} is $G$-equivariant.
\end{lemma}
\begin{proof}
	For each $\sigma\in\Sigma$, let $\mathfrak{a}_{\sigma}\subset {K}[S_{\sigma}]$ denote the ideal cutting out $\imath(X)\cap Y_{\sigma}$. A point $v$ in $\text{Trop}(X,\imath)\cap\text{Trop}(Y_{\sigma})$ is given by a semigroup homomorphism
	$$
	v:S_{\sigma}\longrightarrow\RR\sqcup\{\infty\}
	$$
that admits at least one factorization
	$$
	S_{\sigma}\mono{K}[S_{\sigma}]\xrightarrow{\ |-|_{y}\ }\RR_{\ge0}\xrightarrow{-\log\ \ }\RR\sqcup\{\infty\},
	$$
for some $y\in Y^{\an}_{\sigma}$, such that $|\mathfrak{a}_{\sigma}|_{y}=\{0\}$. The assumption that $\imath$ is $G$-equivariant implies that the isomorphism \eqref{equation: action of g on semigroup algebras} satisfies $g(\mathfrak{a}_{g^{\ast\!}(\sigma)})=\mathfrak{a}_{\sigma}$. Thus if $K_{0}$ is Henselian, commutativity of the diagrams \eqref{equation: commuting action on semigroup algebras} and \eqref{equation: commuting needle point} implies that $g^{\ast}(v)$ lies inside $\text{Trop}({X},\imath)\cap\text{Trop}(Y_{g^{\ast\!}(\sigma)})$. That \eqref{equation: tropicalization map on X} is $G$-equivariant follows from commutativity of the diagram \eqref{equation: commuting needle point}. 
\end{proof}

\begin{definition}
\normalfont
	If $\imath:{X}\mono Y_{\Sigma}$ is a $G$-equivariant closed embedding into a $G^{\op}$-twisted toric ${K}$-variety, then the {\em $G$-equivariant tropicalization of ${X}$} ({\em with respect to $\imath$}) is the extended tropicalization $\text{Trop}({X},\imath)$ equipped with the $G^{\op}$-action provided by Lemma \eqref{lemma: invariance and equivariance}. We denote it
	$$
	\text{Trop}_{G}({X},\imath),
	$$
and refer to the map \eqref{equation: tropicalization map on X} as its {\em $G^{\op}$-equivariant tropicalizaiton map}.
\end{definition}

\begin{example}\label{Galois equi trop of P^1}
\normalfont
{\bf A Galois-equivariant tropicalization of $\pmb{\PP^{1}}$.}
	Let $K_{0}=\CC(\!(t)\!)$ with absolute Galois group $G=\text{Gal}\big(\CC(\!(t^{\QQ})\!)\big/\CC(\!(t)\!)\big)$ as in Example \ref{example: Brauer-Severi A}. Let $X_{0}$ be the projective line over $\CC(\!(t)\!)$:
	$$
	X_{0}
	\ \Def\ 
	\text{Proj}_{\ \!}\CC(\!(t)\!)[x_{0},x].
	$$
Consider the surjective homomorphism $G\epi\text{Gal}\big(\CC(\!(t^{\frac{1}{2}})\!)\big/\CC(\!(t)\!)\big)\cong\ZZ/2\ZZ$. Define $M\Def\ZZ u\oplus\ZZ v\oplus\ZZ w$ with $G$-action factoring trough the $\ZZ/2\ZZ$-action that exchanges $u$ and $v$. Note that the plane in $N_{\RR}$ dual to $u-v\in M$ is the maximal $\ZZ/2\ZZ$-fixed subspace of $N_{\RR}$. Let $\Sigma$ be the $G^{\op}$-invariant fan in $N_{\RR}$ with rays spanned by the vectors $e_{1}$, $e_{2}$, $e_{3}$, and $-(e_{1}+e_{2}+e_{3})$. It determines the toric $\CC(\!(t^{\QQ})\!)$-variety $Y_{\Sigma}\cong\PP^{3}_{\sm{\CC(\!(t^{\QQ})\!)}}$ with twisted $G^{\op}$-action. Give ${X}=X_{0,\sm{\CC(\!(t^{\QQ})\!)}}$ the canonical $G$-action. Then we have a $G$-equivariant closed embedding
	\begin{equation}\label{equation: closed embedding of P^1}
	\begin{array}{rcl}
	\imath\ \!:\ {X}
	&
	\!\!\!\!\mono\!\!\!\!\!\!
	&
	Y_{\Sigma}
	\\[0pt]
	\mbox{{\smaller\smaller $x-t^{\frac{1}{2}}$}}\!\!\!\!\!\!
	&
	\reflectbox{{\smaller\smaller $\longmapsto$}}\!\!\!\!\!
	&
	\mbox{{\smaller\smaller $\chi^{u}$}}
	\\[-2pt]
	\mbox{{\smaller\smaller $x+t^{\frac{1}{2}}$}}\!\!\!\!\!\!
	&
	\reflectbox{{\smaller\smaller $\longmapsto$}}\!\!\!\!\!
	&
	\mbox{{\smaller\smaller $\chi^{v}$}}
	\\[-2pt]
	\mbox{{\smaller\smaller $x-1$}}\!\!\!\!\!\!
	&
	\reflectbox{{\smaller\smaller $\longmapsto$}}\!\!\!\!\!
	&
	\mbox{{\smaller\smaller $\chi^{w}$}},
	\end{array}
	\end{equation}
hence a $G$-equivariant tropicalization $\text{trop}:{X}^{\an}\longrightarrow\text{Trop}_{G}({X},\imath)$. See Figure \ref{figure: trop with action pictured} for a depiction.
	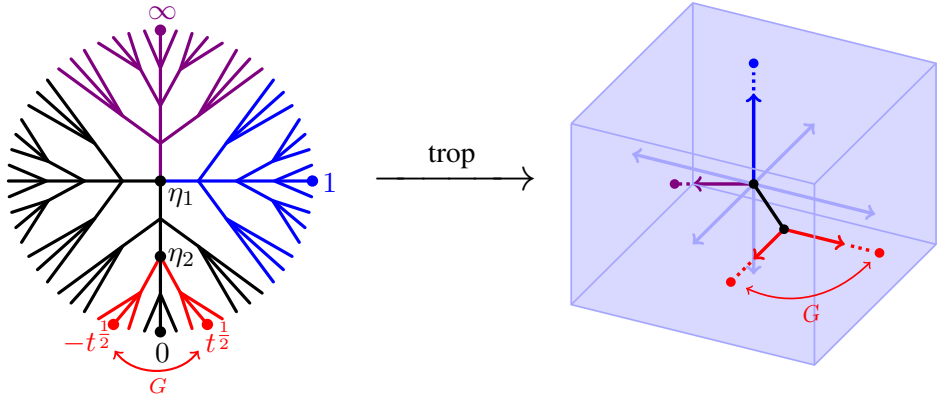
\begin{figure}[h!]
	$$
	\begin{array}{ccc}
	\begin{xy}
	(0,-2.5)*+{
	\begin{tikzpicture}
	\draw[black, very thick, line cap=round, rounded corners=.1mm] (0,0) -- (0,-2);
	\draw[violet, very thick, line cap=round, rounded corners=.1mm] (0,0) -- (0,2);
	\draw[black, very thick, line cap=round, rounded corners=.1mm] (-2,0) -- (0,0);
	\draw[blue, very thick, line cap=round, rounded corners=.1mm] (0,0) -- (2,0);
	\draw[blue, very thick, line cap=round, rounded corners=.1mm] ({.5*cos(0)},{.5*sin(0)}) -- ({cos(30)},{sin(30)}) -- ({2*cos(30)},{2*sin(30)});
	\draw[blue, very thick, line cap=round, rounded corners=.1mm] ({2*cos(42)},{2*sin(42)}) -- ({cos(30)},{sin(30)}) -- ({2*cos(36)},{2*sin(36)});
	\draw[blue, very thick, line cap=round, rounded corners=.1mm] ({.5*cos(0)},{.5*sin(0)}) -- ({cos(-30)},{sin(-30)}) -- ({2*cos(-30)},{2*sin(-30)});
	\draw[blue, very thick, line cap=round, rounded corners=.1mm] ({2*cos(-42)},{2*sin(-42)}) -- ({cos(-30)},{sin(-30)}) -- ({2*cos(-36)},{2*sin(-36)});
	\draw[blue, very thick, line cap=round, rounded corners=.1mm] ({cos(0)},{sin(0)}) -- ({1.5*cos(10)},{1.5*sin(10)}) -- ({2*cos(12)},{2*sin(12)});
	\draw[blue, very thick, line cap=round, rounded corners=.1mm] ({2*cos(24)},{2*sin(24)}) -- ({1.5*cos(10)},{1.5*sin(10)}) -- ({2*cos(18)},{2*sin(18)});
	\draw[blue, very thick, line cap=round, rounded corners=.1mm] ({cos(0)},{sin(0)}) -- ({1.5*cos(-10)},{1.5*sin(-10)}) -- ({2*cos(-12)},{2*sin(-12)});
	\draw[blue, very thick, line cap=round, rounded corners=.1mm] ({2*cos(-24)},{2*sin(-24)}) -- ({1.5*cos(-10)},{1.5*sin(-10)}) -- ({2*cos(-18)},{2*sin(-18)});
	\draw[blue, very thick, line cap=round, rounded corners=.1mm] ({2*cos(6)},{2*sin(6)}) -- ({1.5*cos(0)},{1.5*sin(0)}) -- ({2*cos(-6)},{2*sin(-6)});
	\draw[black, very thick, line cap=round, rounded corners=.1mm] ({.5*cos(0-90)},{.5*sin(0-90)}) -- ({cos(30-90)},{sin(30-90)}) -- ({2*cos(30-90)},{2*sin(30-90)});
	\draw[black, very thick, line cap=round, rounded corners=.1mm] ({2*cos(42-90)},{2*sin(42-90)}) -- ({cos(30-90)},{sin(30-90)}) -- ({2*cos(36-90)},{2*sin(36-90)});
	\draw[black, very thick, line cap=round, rounded corners=.1mm] ({.5*cos(0-90)},{.5*sin(0-90)}) -- ({cos(-30-90)},{sin(-30-90)}) -- ({2*cos(-30-90)},{2*sin(-30-90)});
	\draw[black, very thick, line cap=round, rounded corners=.1mm] ({2*cos(-42-90)},{2*sin(-42-90)}) -- ({cos(-30-90)},{sin(-30-90)}) -- ({2*cos(-36-90)},{2*sin(-36-90)});
	\draw[red, very thick, line cap=round, rounded corners=.1mm] ({cos(0-90)},{sin(0-90)}) -- ({1.5*cos(10-90)},{1.5*sin(10-90)}) -- ({2*cos(12-90)},{2*sin(12-90)});
	\draw[red, very thick, line cap=round, rounded corners=.1mm] ({2*cos(24-90)},{2*sin(24-90)}) -- ({1.5*cos(10-90)},{1.5*sin(10-90)}) -- ({2*cos(18-90)},{2*sin(18-90)});
	\draw[red, very thick, line cap=round, rounded corners=.1mm] ({cos(0-90)},{sin(0-90)}) -- ({1.5*cos(-10-90)},{1.5*sin(-10-90)}) -- ({2*cos(-12-90)},{2*sin(-12-90)});
	\draw[red, very thick, line cap=round, rounded corners=.1mm] ({2*cos(-24-90)},{2*sin(-24-90)}) -- ({1.5*cos(-10-90)},{1.5*sin(-10-90)}) -- ({2*cos(-18-90)},{2*sin(-18-90)});
	\draw[black, very thick, line cap=round, rounded corners=.1mm] ({2*cos(6-90)},{2*sin(6-90)}) -- ({1.5*cos(0-90)},{1.5*sin(0-90)}) -- ({2*cos(-6-90)},{2*sin(-6-90)});
	\draw[black, very thick, line cap=round, rounded corners=.1mm] ({.5*cos(0+180)},{.5*sin(0+180)}) -- ({cos(30+180)},{sin(30+180)}) -- ({2*cos(30+180)},{2*sin(30+180)});
	\draw[black, very thick, line cap=round, rounded corners=.1mm] ({2*cos(42+180)},{2*sin(42+180)}) -- ({cos(30+180)},{sin(30+180)}) -- ({2*cos(36+180)},{2*sin(36+180)});
	\draw[black, very thick, line cap=round, rounded corners=.1mm] ({.5*cos(0+180)},{.5*sin(0+180)}) -- ({cos(-30+180)},{sin(-30+180)}) -- ({2*cos(-30+180)},{2*sin(-30+180)});
	\draw[black, very thick, line cap=round, rounded corners=.1mm] ({2*cos(-42+180)},{2*sin(-42+180)}) -- ({cos(-30+180)},{sin(-30+180)}) -- ({2*cos(-36+180)},{2*sin(-36+180)});
	\draw[black, very thick, line cap=round, rounded corners=.1mm] ({cos(0+180)},{sin(0+180)}) -- ({1.5*cos(10+180)},{1.5*sin(10+180)}) -- ({2*cos(12+180)},{2*sin(12+180)});
	\draw[black, very thick, line cap=round, rounded corners=.1mm] ({2*cos(24+180)},{2*sin(24+180)}) -- ({1.5*cos(10+180)},{1.5*sin(10+180)}) -- ({2*cos(18+180)},{2*sin(18+180)});
	\draw[black, very thick, line cap=round, rounded corners=.1mm] ({cos(0+180)},{sin(0+180)}) -- ({1.5*cos(-10+180)},{1.5*sin(-10+180)}) -- ({2*cos(-12+180)},{2*sin(-12+180)});
	\draw[black, very thick, line cap=round, rounded corners=.1mm] ({2*cos(-24+180)},{2*sin(-24+180)}) -- ({1.5*cos(-10+180)},{1.5*sin(-10+180)}) -- ({2*cos(-18+180)},{2*sin(-18+180)});
	\draw[black, very thick, line cap=round, rounded corners=.1mm] ({2*cos(6+180)},{2*sin(6+180)}) -- ({1.5*cos(0+180)},{1.5*sin(0+180)}) -- ({2*cos(-6+180)},{2*sin(-6+180)});
	\draw[violet, very thick, line cap=round, rounded corners=.1mm] ({.5*cos(0+90)},{.5*sin(0+90)}) -- ({cos(30+90)},{sin(30+90)}) -- ({2*cos(30+90)},{2*sin(30+90)});
	\draw[violet, very thick, line cap=round, rounded corners=.1mm] ({2*cos(42+90)},{2*sin(42+90)}) -- ({cos(30+90)},{sin(30+90)}) -- ({2*cos(36+90)},{2*sin(36+90)});
	\draw[violet, very thick, line cap=round, rounded corners=.1mm] ({.5*cos(0+90)},{.5*sin(0+90)}) -- ({cos(-30+90)},{sin(-30+90)}) -- ({2*cos(-30+90)},{2*sin(-30+90)});
	\draw[violet, very thick, line cap=round, rounded corners=.1mm] ({2*cos(-42+90)},{2*sin(-42+90)}) -- ({cos(-30+90)},{sin(-30+90)}) -- ({2*cos(-36+90)},{2*sin(-36+90)});
	\draw[violet, very thick, line cap=round, rounded corners=.1mm] ({cos(0+90)},{sin(0+90)}) -- ({1.5*cos(10+90)},{1.5*sin(10+90)}) -- ({2*cos(12+90)},{2*sin(12+90)});
	\draw[violet, very thick, line cap=round, rounded corners=.1mm] ({2*cos(24+90)},{2*sin(24+90)}) -- ({1.5*cos(10+90)},{1.5*sin(10+90)}) -- ({2*cos(18+90)},{2*sin(18+90)});
	\draw[violet, very thick, line cap=round, rounded corners=.1mm] ({cos(0+90)},{sin(0+90)}) -- ({1.5*cos(-10+90)},{1.5*sin(-10+90)}) -- ({2*cos(-12+90)},{2*sin(-12+90)});
	\draw[violet, very thick, line cap=round, rounded corners=.1mm] ({2*cos(-24+90)},{2*sin(-24+90)}) -- ({1.5*cos(-10+90)},{1.5*sin(-10+90)}) -- ({2*cos(-18+90)},{2*sin(-18+90)});
	\draw[violet, very thick, line cap=round, rounded corners=.1mm] ({2*cos(6+90)},{2*sin(6+90)}) -- ({1.5*cos(0+90)},{1.5*sin(0+90)}) -- ({2*cos(-6+90)},{2*sin(-6+90)});
	\draw[red, thick, <->] (-.6,-2.25) to [out=-55, in=-125] (.5,-2.25);
	\fill[blue] (2,0) circle (.075);
	\fill[black] (0,0) circle (.075);
	\fill[black] (0,-1) circle (.075);
	\fill[black] (0,-2) circle (.075);
	\fill[violet] (0,2) circle (.075);
	\fill[red] ({2*cos(18-90)},{2*sin(18-90)}) circle (.075);
	\fill[red] ({2*cos(-18-90)},{2*sin(-18-90)}) circle (.075);
	\end{tikzpicture}
	};
	(2.5,-2)*+{\eta_{1}};
	(2.5,-10)*+{\eta_{2}};
	(22,0)*+{{\color{blue} 1}};
	(0,-22.5)*+{0};
	(0,22.5)*+{{\color{violet} \infty}};
	(7.5,-20.5)*+{{\color{red} t^{\frac{1}{2}}}};
	(-9.75,-20.5)*+{{\color{red} -t^{\frac{1}{2}}}};
	(0,-27)*+{{\color{red} {}_{G\ }}};
	\end{xy}
	&
	\mbox{{\larger\larger\larger $\xrightarrow{\ \ \ \ \mbox{{\smaller\smaller\smaller trop}}\ \ \ \ }$}}
	&
	\scalebox{.8}{$
	\begin{xy}
	(0,0)*+{
	\begin{tikzpicture}
	\fill[blue!15] (-2,1) -- (-2,-2) -- (2,-3) -- (4,-1) -- (4,2) -- (0,3);
	\draw[violet, ultra thick, ->] (1,0) -- (0,0);
	\draw[violet, ultra thick, dotted] (1-1.05,0) -- (1-1.25,0);
	\draw[blue!35, thick] (-2,-2) -- (0,0) -- (0,3);
	\draw[blue!35, thick] (0,0) -- (4,-1);
	\draw[blue!35, ultra thick, <->] (0,-1) -- (2,1);
	\draw[blue!35, ultra thick, <->] (1,1.5) -- (1,-1.5);
	\draw[blue!35, ultra thick, <->] (-1,.5) -- (3,-.5);
	\draw[black, ultra thick] (1,0) -- (1.5,-.75);
	\draw[blue, ultra thick, ->] (1,0) -- (1,1.5);
	\draw[blue, ultra thick, dotted] (1,1.55) -- (1,1.85);
	\draw[red, ultra thick, <->] (1,-1.25) -- (1.5,-.75) -- (2.5,-1);
	\draw[red, ultra thick, dotted] (1.5-1.1*.5,-.75-1.1*.5) -- (1.5-1.55*.5,-.75-1.55*.5);
	\draw[red, ultra thick, dotted] (1.5+1.1*1,-.75-1.1*.25) -- (1.5+1.4*1,-.75-1.4*.25);
	\draw[blue!35, thick] (-2,1) -- (-2,-2) -- (2,-3) -- (4,-1) -- (4,2) -- (0,3) -- (-2,1);
	\draw[blue!35, thick] (-2,1) -- (2,0) -- (2,-3);
	\draw[blue!35, thick] (2,0) -- (4,2);
	\fill[blue] (1,2) circle (.08);
	\fill[red] (1.5-1.75*.5,-.75-1.75*.5) circle (.08);
	\fill[red] (1.5+1.57*1,-.75-1.57*.25) circle (.08);
	\fill[violet] (1-1.3,0) circle (.08);
	\fill[black] (1,0) circle (.08);
	\fill[black] (1.5,-.75) circle (.08);
	\draw[red, thick, <->] (1.5-1.75*.5+.25,-.75-1.75*.5-.1) to [out=-25, in=-135] (1.5+1.57*1-.15,-.75-1.57*.25-.15);
	\end{tikzpicture}
	};
	(9.5,-21.5)*+{\rotatebox{10}{${\color{red} G}$}};
	\end{xy}
	$}
	\end{array}
	$$
\caption{The Berkovich analytification of the projective line over $\CC(\!(t^{\QQ})\!)$ and its tropicalization with respect to the closed embedding \eqref{equation: closed embedding of P^1}. Here $\eta_{1}$ denotes the weight-$1$ Gauss norm and $\eta_{2}$ denotes the weight-$e^{-\frac{1}{2}}$ Gauss norm.}
\label{figure: trop with action pictured}
\end{figure}
	The action $G^{\op}\ \lact\ {X}^{\an}$ is quite complex, but we can describe pieces of it. Let $B_{r}\big(\pm t^{\frac{1}{2}}\big)^{\circ}$ denote the open analytic disks of radius $r=e^{-\frac{1}{2}}$ in ${X}^{\an}$, centered at $\pm t^{\frac{1}{2}}$ (red in Figure \ref{figure: trop with action pictured}). Then
	$$
	B_{r}\big(\!\pm t^{\frac{1}{2}}\big)^{\!\circ}\cap\ \CC(\!(t)\!)^{\text{sep}}
	\ \ =\ \ 
	\Big\{
	\pm t^{\frac{1}{2}}\!+\!\!\!\sum_{\sm{q\!>\!1/2}}\!\!c_{q}\ \!t^{q}\in\CC(\!(t)\!)^{\text{sep}}
	\ \!\Big\}.
	$$
This implies that every element $g\in G$ mapping to $1$ under $G\epi\ZZ/2\ZZ$ induces a homeomorphism
	$$
	g^{\ast}:B_{r}\big(t^{\frac{1}{2}}\big)^{\!\circ}\xrightarrow{\ \ \sim\ \ }B_{r}\big(\!-t^{\frac{1}{2}}\big)^{\!\circ}.
	$$
The open analytic disk $B_{1}(1)^{\circ}$ (blue in Figure \ref{figure: trop with action pictured}) and the open analytic domain ${X}^{\an}\backslash B_{1}(1)$ (purple in Figure \ref{figure: trop with action pictured}) are each invariant under the full Galois group $G$. All points of ${X}^{\an}$ lying in the dense torus, except those in the open disks $B_{r}\big(\pm t^{\frac{1}{2}}\big)^{\!\circ}$, map into the $\ZZ/2\ZZ$-fixed plane in $N_{\RR}$. The $G$-action on $N_{\RR}$ exchanges the images of the open disks $B_{r}\big(\pm t^{\frac{1}{2}}\big)^{\!\circ}$.
\end{example}

\begin{corollary}\label{corollary: equivariance of tropicalization map}
\normalfont
	Assume that $K_{0}$ is Henselian and equip $X^{\an}$ with its canonical $G^{\op}$-action. Then for any system $\mathcal{S}_{G}$ of $G$-equivariant toric embeddings, the continuous map \eqref{equation: canonical G-equivariant map} is $G$-equivariant.
\hfill
$\square$
\end{corollary}

\end{subsection}


\begin{subsection}{The polyhedral structure on a Galois equivariant tropicalization}\label{unramified Galois descended initial degenerations}
	For the remainder of \S\ref{unramified Galois descended initial degenerations}, assume that $K_{0}$ is Henselian. Denote the valuation ring, maximal ideal, and residue field of ${K}=K^{\text{sep}}_{0}$
	$$
	\begin{array}{rcl}
	{K}^{\circ}
	&
	\!\!\!=\!\!\!
	& 
	\big\{
	a\in {K}:|a|\le 1
	\big\},
	\\[6pt]
	{K}^{\circ\circ}
	&
	\!\!\!=\!\!\!
	& 
	\big\{
	a\in {K}:|a|< 1
	\big\},
	\\[6pt]
	\widetilde{K}
	&
	\!\!\!=\!\!\!
	&
	{K}^{\circ}\!\big/{K}^{\circ\circ}.
	\end{array}
	$$
Suppose given a $K_{0}$-variety $X_{0}$ and a closed embedding $\imath:{X}\mono Y_{\Sigma}$ of $X=X_{0,{K}}$ into a toric ${K}$-variety $Y_{\Sigma}$. Equip $Y_{\Sigma}$ with a twisted $G^{\op}$-action $G^{\op}\ \lact\ Y_{\Sigma}$.
	
	As in the proof of Lemma \ref{lemma: invariance and equivariance}, for each $\sigma\in \Sigma$, let $\mathfrak{a}_{\sigma}\subset {K}[M]$ denote the ideal cutting out $\imath(X)\cap Y_{\sigma}$. Each point $v\in N_{\RR}$ determines a $K^{\circ}$-model $\mathscr{X}^{v}$ of $X$ glued from the spectra of the $K^{\circ}$-algebras ${K}[S_{\sigma}]^{v}\big/\mathfrak{a}^{v}_{\sigma}$, where
	$$
	{K}[S_{\sigma}]^{v}
	\ \ \Def\ \ 
	 \big\{a\chi^{u}\in{K}[S_{\sigma}]:\langle u,v\rangle-\log|a|\ge0
	 \big\}
	$$
\begin{center}
and
\end{center}
	$$
	\mathfrak{a}^{v}_{\sigma}
	\ \ \Def\ \ 
	\mathfrak{a}_{\sigma}\cap{K}[S_{\sigma}]^{v}.
	$$
Because each ${K}[S_{\sigma}]^{v}$ is a subring of ${K}[S_{\sigma}]$, it makes sense to say that the special fibers $\mathscr{X}^{v_{1}}_{\sm{\widetilde{K}}}$ and $\mathscr{X}^{v_{2}}_{\sm{\widetilde{K}}}$ are ``equal" or ``not equal" for distinct $v_{1},v_{2}\in N_{\RR}$. For any closed embedding $\imath:X\mono\PP^{n}_{\!\!\sm{K}}$, where $\PP^{n}_{\!\!\sm{K}}$ has its standard toric structure, the {\em Gr\"obner complex} on $\text{Trop}(Y_{\Sigma})$ is the closure of the complex in $N_{\RR}$ whose cells consist of vectors $v$ with fixed $\mathscr{X}^{v}$ \cite[\S2.5]{MS:15} \cite[proof of Theorem 10.14]{Gub:13}.
	
	Define a {\em $G^{\op}$-twisted $\PP^{n}_{\!\!\sm{K}}$} to be any $G^{\op}$-twisted toric ${K}$-variety with underlying ${K}$-variety $Y_{\Sigma}\cong\PP^{n}_{\!\!\sm{K}}$.

\begin{lemma}\label{lemma: invariance of full Groebner complex}
\normalfont
	If $K_{0}$ is Henselian and $\imath:{X}\mono Y_{\Sigma}$ is a $G$-equivariant closed embedding into a $G^{\op}$-twisted $\PP^{n}_{\!\!\sm{K}}$, the Gr\"obner complex in $\text{Trop}(Y_{\Sigma})$ induced by $\imath$ is invariant under the $G^{\op}$-action on $\text{Trop}(Y_{\Sigma})$.
\end{lemma}
\begin{proof}
	If $Y_{\Sigma}$ is a $G^{\op}$-twisted toric variety and the closed embedding $\imath:{X}\mono Y_{\Sigma}$ is $G$-equivariant, then for each $\sigma\in\Sigma$, each $u\in S_{\sigma}$, each $v\in\text{Trop}(Y_{\Sigma})$, and each $g\in G$, invariance of the absolute value on ${K}$ under $g$ implies that
	$$
	\big\langle g(u),v\big\rangle-\log\big|g(a)\big|\ge0
	\ \ \ \ \ \ \mbox{if and only if}\ \ \ \ \ \ 
	\big\langle u,g^{\ast}(v)\big\rangle-\log|a|\ge0.
	$$
Thus the automorphism $g\ \lact\ {K}[M]$ restricts to an isomorphism
	$$
	g:{K}[S_{g^{\ast\!}(\sigma)}]^{g^{\ast\!}(v)}\xrightarrow{\ \sim\ }{K}[S_{\sigma}]^{v},
	$$
and for any pair of points $v_{1},v_{2}\in N_{\RR}$, we have ${K}[S_{\sigma}]^{v_{1}}={K}[S_{\sigma}]^{v_{2}}$ inside ${K}[S_{\sigma}]$ if and only if ${K}[S_{g^{\ast\!}(\sigma)}]^{g^{\ast\!}(v_{1})}={K}[S_{g^{\ast\!}(\sigma)}]^{g^{\ast\!}(v_{2})}$ inside ${K}[S_{g^{\ast\!}(\sigma)}]$. Because $g(\mathfrak{a}_{g^{\ast\!}(\sigma)})=\mathfrak{a}_{\sigma}$, this implies that $v_{1}$ and $v_{2}$ lie in the same Gr\"obner cell in $N_{\RR}$ if an only if $g^{\ast\!}(v_{1})$ and $g^{\ast\!}(v_{2})$ lie in the same Gr\"obner cell in $N_{\RR}$.
\end{proof}

\begin{corollary}\label{corollary: existence of G-invariant polyhedral complex}
\normalfont
	If $K_{0}$ is Henselian and $\imath:{X}\mono Y_{\Sigma}$ is a $G$-equivariant closed embedding into a $G^{\op}$-twisted $\PP^{n}_{\!\!\sm{K}}$, then $\text{Trop}_{G}({X},\imath)$ admits a $G^{\op}$-invariant polyhedral decomposition.
\end{corollary}
\begin{proof}
	This follows immediately from Lemma \ref{lemma: invariance of full Groebner complex} and \cite[Theorem 10.14]{Gub:13}.
\end{proof}

\begin{corollary}\label{corollary: existence of fine enough G-invariant polyhedral complex}
\normalfont
	If $K_{0}$ is Henselian and $\imath:{X}\mono Y_{\Sigma}$ is a $G$-equivariant closed embedding into a $G^{\op}$-twisted $\PP^{n}_{\!\!\sm{K}}$, then there exists a $G$-invariant polyhedral decomposition of $\text{Trop}_{G}(X,\imath)$ that refines the polyhedral decomposition $\{\overline{\ \!\sigma\ \!}\}_{\sigma\in\Sigma}$.
\end{corollary}
\begin{proof}
	Let $\mathcal{P}$ be the Gr\"obner complex on $\text{Trop}({X},\imath)$. Then the set
	$
	\mathcal{Q}
	\Def
	\big\{\ \!
	P\cap\overline{\ \!\sigma\ \!}
	\ \!:\ \!
	P\in\mathcal{P}\ \text{and}\ \overline{\ \!\sigma\ \!}\in\Sigma
	\ \!\big\}
	$
is a polyhedral decomposition with the required property.
\end{proof}

\begin{proof}[{\it Proof of Theorem \ref{theorem: main D}}]
	By Corollary \ref{corollary: existence of G-invariant polyhedral complex}, $\text{Trop}({X},\imath)$ admits a polyhedral decomposition that is $G^{\op}$-invariant. For each polyhedral cell $\Delta$ in this polyhedral decomposition and for each integer $p\ge0$, every element $g\in G$ induces isomorphisms
	\begin{equation}\label{equation: action on cochain complex terms}
	g:\mathcal{F}^{p}(\Delta)\xrightarrow{\ \sim\ }\mathcal{F}^{p}(g^{\ast}\!\Delta)
	\end{equation}
between the terms $\mathcal{F}^{p}(\Delta)$ and $\mathcal{F}(g^{\ast}\!\Delta)$ in the cochain group $C^{\text{dim}_{\ \!}\Delta}\big(\text{Trop}({X},\imath),\mathcal{F}^{p}\big)$ as defined in \cite[Definition 16 on]{IKMZ:16}. Thus for each integer $q\ge0$, the cochain group $C^{q}\big(\text{Trop}({X},\imath),\mathcal{F}^{p}\big)$ picks up a $G$-action. This action \eqref{equation: action on cochain complex terms} commutes with the restriction morphisms ``$\iota^{\ast}$" defined in \cite[Equation (8)]{IKMZ:16}. Thus the coboundary morphisms \cite[\S2.3, after Equation (8)]{IKMZ:16} on $C^{q}\big(\text{Trop}({X},\imath),\mathcal{F}^{p}\big)$ are $G$-equivariant, and we obtain a $G$-action on each tropical cellular cohomology group:
	$$
	G\ \ \lact\ \ H^{p,q}_{\text{trop}}\big(\ \!\text{Trop}({X},\imath),F\ \!\big)
	\ \ \ \mbox{for each}\ \ \ p,q\ge0.
	$$
The representation on tropical homology arises in the same way after dualizing $\mathcal{F}^{p}(\Delta)\mapsto\mathcal{F}_{p}(\Delta)$.
\end{proof}

\begin{example}\label{example: Brauer-Severi C}
\normalfont
{\bf A Galois representation in tropical cellular cohomology.}
	Let $K_{0}=\CC(\!(t)\!)$ with separable closure $K=\CC(\!(t^{\QQ})\!)$ and Galois group $G=\text{Gal}\big(\CC(\!(t^{\QQ})\!)\big/\CC(\!(t)\!)\big)$. Let $Y_{0}$ be the Brauer-Severi $\CC(\!(t)\!)$-variety constructed in Example \ref{example: Brauer-Severi A}. If we give $M=\text{Hom}_{\ZZ}(N,\ZZ)$ the basis dual to the standard basis on $N=\ZZ^{2}$, then the action of the matrix $\left(\begin{smallmatrix}0 & -1 \\ 1 & -1\end{smallmatrix}\right)$ on $N$ in Example \ref{example: Brauer-Severi A} is dual to the $3$-periodic action of its transpose $\left(\begin{smallmatrix}\ 0\! & \ 1\! \\ \!-1 & \!-1\end{smallmatrix}\right)$ on $M$. This action induces a twisted $\ZZ/3\ZZ$-action, and thus a $G$-action, on $\CC(\!(t^{\frac{1}{3}})\!)[M]=\CC(\!(t^{\frac{1}{3}})\!)[x^{\pm1},y^{\pm1}]$. Consider the Laurent polynomial
	\begin{equation}\label{equation: genus 3 Laurent polynomial}
	\scalebox{.9}{$
	\begin{array}{c}
	t^{12}\ \!\big(x^{-2}y^{-2}+x^{-2}y^{4}+x^{4}y^{-2})+t^{7}\ \!\big(x^{-1}y^{-2}+x^{-2}y^{3}+x^{3}y^{-1}\big)+t^{7}\ \!\big(x^{3}y^{-2}+x^{-2}y^{-1}+x^{-1}y^{3}\big)
	\\[7pt]
	+\ t^{4}\ \!\big(y^{-2}+x^{-2}y^{2}+x^{2}\big)+t^{4}\ \!\big(y^{2}+x^{2}y^{-2}+x^{-2}\big)+t^{3}\ \!\big(xy^{-2}+x^{-2}y+xy\big)
	\\[7pt]
	+\ t^{3}\ \!\big(x^{-1}y^{2}+x^{2}y^{-1}+x^{-1}y^{-1}\big)+t\ \!\big(y^{-1}+x^{-1}y+x\big)+t\ \!\big(y+xy^{-1}+x^{-1}\big)+1
	\end{array}
	$}
	\end{equation}
in $\CC(\!(t^{\frac{1}{3}})\!)[x^{\pm1},y^{\pm1}]$. It cuts out a smooth, degree-$6$, genus-$10$ curve $X$ inside $Y_{\Sigma}\cong\PP^{2}_{\sm{\CC(\!(t^{\QQ})\!)}}$. Each of the nine trinomials in parentheses in \eqref{equation: genus 3 Laurent polynomial} constitutes an orbit under the $\ZZ/3\ZZ$-action on $\CC(\!(t^{\frac{1}{3}})\!)[x^{\pm1},y^{\pm1}]$, making $X$ a $G^{\op}$-invariant subvariety of the twisted toric $\CC(\!(t^{\QQ})\!)$-variety $G^{\op}\ \lact\ Y_{\Sigma}$. It arises as the pullback of a $\CC(\!(t)\!)$-curve $X_{0}$ inside the resulting Brauer-Severi $\CC(\!(t)\!)$-variety $Y_{0}$.
	
	The $G$-equivariant tropicalization of ${X}$ is pictured in Figure \ref{figure: example of tropicalization} below. Note the apparent $\ZZ/3\ZZ$-symmetry.
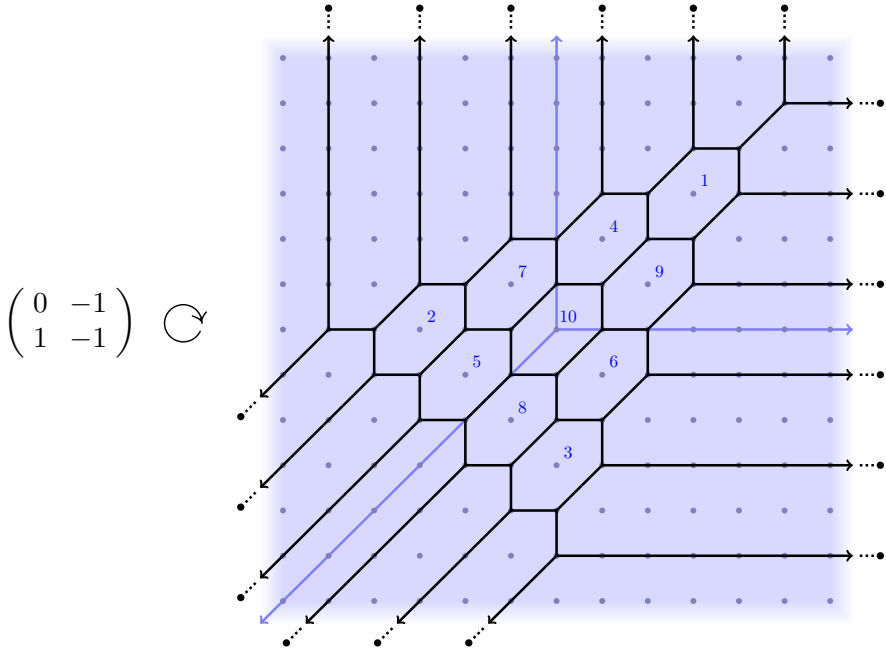
\begin{figure}[h!]
	$$
	\!\!\!\!\!\!\!\!\!\!\!\!\!\!\!\!\!\!\!\!\!\!\!\!\!\!\!\!\!\!\!\!\!\!\!\!\!\!\!\!\!\!\!
	\left(\!
	\begin{array}{cc}
	0 & \!-1 \\
	1 & \!-1
	\end{array}
	\!\right)
	\ \ 
	\raisebox{-5pt}{
	\scalebox{2}{$\lact$}
	}
	\scalebox{.6}{$
	\ \ \ \ \ 
	\begin{xy}
	(0,0)*+{
	\begin{tikzpicture}
	\fill[blue!15, path fading=north] (-6,6) -- (6,6) -- (6.5,6.5) -- (-6.5,6.5);
	\fill[blue!15, path fading=south] (-6,-6) -- (6,-6) -- (6.5,-6.5) -- (-6.5,-6.5);
	\fill[blue!15, path fading=east] (6,6) -- (6,-6) -- (6.5,-6.5) -- (6.5,6.5);
	\fill[blue!15, path fading=west] (-6,6) -- (-6,-6) -- (-6.5,-6.5) -- (-6.5,6.5);
	\fill[blue!15] (-6,-6) -- (6,-6) -- (6,6) -- (-6,6);
	\draw[blue!50, ultra thick, ->] (0,0) -- (6.5,0);
	\draw[blue!50, ultra thick, ->] (0,0) -- (0,6.5);
	\draw[blue!50, ultra thick, ->] (0,0) -- (-6.5,-6.5);
	\fill[black!50!blue!50] (0,0) circle (.065);
	\fill[black!50!blue!50] (1,0) circle (.065);
	\fill[black!50!blue!50] (2,0) circle (.065);
	\fill[black!50!blue!50] (3,0) circle (.065);
	\fill[black!50!blue!50] (4,0) circle (.065);
	\fill[black!50!blue!50] (5,0) circle (.065);
	\fill[black!50!blue!50] (6,0) circle (.065);
	\fill[black!50!blue!50] (-1,0) circle (.065);
	\fill[black!50!blue!50] (-2,0) circle (.065);
	\fill[black!50!blue!50] (-3,0) circle (.065);
	\fill[black!50!blue!50] (-4,0) circle (.065);
	\fill[black!50!blue!50] (-5,0) circle (.065);
	\fill[black!50!blue!50] (-6,0) circle (.065);
	\fill[black!50!blue!50] (0,1) circle (.065);
	\fill[black!50!blue!50] (1,1) circle (.065);
	\fill[black!50!blue!50] (2,1) circle (.065);
	\fill[black!50!blue!50] (3,1) circle (.065);
	\fill[black!50!blue!50] (4,1) circle (.065);
	\fill[black!50!blue!50] (5,1) circle (.065);
	\fill[black!50!blue!50] (6,1) circle (.065);
	\fill[black!50!blue!50] (-1,1) circle (.065);
	\fill[black!50!blue!50] (-2,1) circle (.065);
	\fill[black!50!blue!50] (-3,1) circle (.065);
	\fill[black!50!blue!50] (-4,1) circle (.065);
	\fill[black!50!blue!50] (-5,1) circle (.065);
	\fill[black!50!blue!50] (-6,1) circle (.065);
	\fill[black!50!blue!50] (0,2) circle (.065);
	\fill[black!50!blue!50] (1,2) circle (.065);
	\fill[black!50!blue!50] (2,2) circle (.065);
	\fill[black!50!blue!50] (3,2) circle (.065);
	\fill[black!50!blue!50] (4,2) circle (.065);
	\fill[black!50!blue!50] (5,2) circle (.065);
	\fill[black!50!blue!50] (6,2) circle (.065);
	\fill[black!50!blue!50] (-1,2) circle (.065);
	\fill[black!50!blue!50] (-2,2) circle (.065);
	\fill[black!50!blue!50] (-3,2) circle (.065);
	\fill[black!50!blue!50] (-4,2) circle (.065);
	\fill[black!50!blue!50] (-5,2) circle (.065);
	\fill[black!50!blue!50] (-6,2) circle (.065);
	\fill[black!50!blue!50] (0,-1) circle (.065);
	\fill[black!50!blue!50] (1,-1) circle (.065);
	\fill[black!50!blue!50] (2,-1) circle (.065);
	\fill[black!50!blue!50] (3,-1) circle (.065);
	\fill[black!50!blue!50] (4,-1) circle (.065);
	\fill[black!50!blue!50] (5,-1) circle (.065);
	\fill[black!50!blue!50] (6,-1) circle (.065);
	\fill[black!50!blue!50] (-1,-1) circle (.065);
	\fill[black!50!blue!50] (-2,-1) circle (.065);
	\fill[black!50!blue!50] (-3,-1) circle (.065);
	\fill[black!50!blue!50] (-4,-1) circle (.065);
	\fill[black!50!blue!50] (-5,-1) circle (.065);
	\fill[black!50!blue!50] (-6,-1) circle (.065);
	\fill[black!50!blue!50] (0,-2) circle (.065);
	\fill[black!50!blue!50] (1,-2) circle (.065);
	\fill[black!50!blue!50] (2,-2) circle (.065);
	\fill[black!50!blue!50] (3,-2) circle (.065);
	\fill[black!50!blue!50] (4,-2) circle (.065);
	\fill[black!50!blue!50] (5,-2) circle (.065);
	\fill[black!50!blue!50] (6,-2) circle (.065);
	\fill[black!50!blue!50] (-1,-2) circle (.065);
	\fill[black!50!blue!50] (-2,-2) circle (.065);
	\fill[black!50!blue!50] (-3,-2) circle (.065);
	\fill[black!50!blue!50] (-4,-2) circle (.065);
	\fill[black!50!blue!50] (-5,-2) circle (.065);
	\fill[black!50!blue!50] (-6,-2) circle (.065);
	\fill[black!50!blue!50] (0,-3) circle (.065);
	\fill[black!50!blue!50] (1,-3) circle (.065);
	\fill[black!50!blue!50] (2,-3) circle (.065);
	\fill[black!50!blue!50] (3,-3) circle (.065);
	\fill[black!50!blue!50] (4,-3) circle (.065);
	\fill[black!50!blue!50] (5,-3) circle (.065);
	\fill[black!50!blue!50] (6,-3) circle (.065);
	\fill[black!50!blue!50] (-1,-3) circle (.065);
	\fill[black!50!blue!50] (-2,-3) circle (.065);
	\fill[black!50!blue!50] (-3,-3) circle (.065);
	\fill[black!50!blue!50] (-4,-3) circle (.065);
	\fill[black!50!blue!50] (-5,-3) circle (.065);
	\fill[black!50!blue!50] (-6,-3) circle (.065);
	\fill[black!50!blue!50] (0,3) circle (.065);
	\fill[black!50!blue!50] (1,3) circle (.065);
	\fill[black!50!blue!50] (2,3) circle (.065);
	\fill[black!50!blue!50] (3,3) circle (.065);
	\fill[black!50!blue!50] (4,3) circle (.065);
	\fill[black!50!blue!50] (5,3) circle (.065);
	\fill[black!50!blue!50] (6,3) circle (.065);
	\fill[black!50!blue!50] (-1,3) circle (.065);
	\fill[black!50!blue!50] (-2,3) circle (.065);
	\fill[black!50!blue!50] (-3,3) circle (.065);
	\fill[black!50!blue!50] (-4,3) circle (.065);
	\fill[black!50!blue!50] (-5,3) circle (.065);
	\fill[black!50!blue!50] (-6,3) circle (.065);
	\fill[black!50!blue!50] (0,4) circle (.065);
	\fill[black!50!blue!50] (1,4) circle (.065);
	\fill[black!50!blue!50] (2,4) circle (.065);
	\fill[black!50!blue!50] (3,4) circle (.065);
	\fill[black!50!blue!50] (4,4) circle (.065);
	\fill[black!50!blue!50] (5,4) circle (.065);
	\fill[black!50!blue!50] (6,4) circle (.065);
	\fill[black!50!blue!50] (-1,4) circle (.065);
	\fill[black!50!blue!50] (-2,4) circle (.065);
	\fill[black!50!blue!50] (-3,4) circle (.065);
	\fill[black!50!blue!50] (-4,4) circle (.065);
	\fill[black!50!blue!50] (-5,4) circle (.065);
	\fill[black!50!blue!50] (-6,4) circle (.065);
	\fill[black!50!blue!50] (0,5) circle (.065);
	\fill[black!50!blue!50] (1,5) circle (.065);
	\fill[black!50!blue!50] (2,5) circle (.065);
	\fill[black!50!blue!50] (3,5) circle (.065);
	\fill[black!50!blue!50] (4,5) circle (.065);
	\fill[black!50!blue!50] (5,5) circle (.065);
	\fill[black!50!blue!50] (6,5) circle (.065);
	\fill[black!50!blue!50] (-1,5) circle (.065);
	\fill[black!50!blue!50] (-2,5) circle (.065);
	\fill[black!50!blue!50] (-3,5) circle (.065);
	\fill[black!50!blue!50] (-4,5) circle (.065);
	\fill[black!50!blue!50] (-5,5) circle (.065);
	\fill[black!50!blue!50] (-6,5) circle (.065);
	\fill[black!50!blue!50] (0,6) circle (.065);
	\fill[black!50!blue!50] (1,6) circle (.065);
	\fill[black!50!blue!50] (2,6) circle (.065);
	\fill[black!50!blue!50] (3,6) circle (.065);
	\fill[black!50!blue!50] (4,6) circle (.065);
	\fill[black!50!blue!50] (5,6) circle (.065);
	\fill[black!50!blue!50] (6,6) circle (.065);
	\fill[black!50!blue!50] (-1,6) circle (.065);
	\fill[black!50!blue!50] (-2,6) circle (.065);
	\fill[black!50!blue!50] (-3,6) circle (.065);
	\fill[black!50!blue!50] (-4,6) circle (.065);
	\fill[black!50!blue!50] (-5,6) circle (.065);
	\fill[black!50!blue!50] (-6,6) circle (.065);
	\fill[black!50!blue!50] (0,-4) circle (.065);
	\fill[black!50!blue!50] (1,-4) circle (.065);
	\fill[black!50!blue!50] (2,-4) circle (.065);
	\fill[black!50!blue!50] (3,-4) circle (.065);
	\fill[black!50!blue!50] (4,-4) circle (.065);
	\fill[black!50!blue!50] (5,-4) circle (.065);
	\fill[black!50!blue!50] (6,-4) circle (.065);
	\fill[black!50!blue!50] (-1,-4) circle (.065);
	\fill[black!50!blue!50] (-2,-4) circle (.065);
	\fill[black!50!blue!50] (-3,-4) circle (.065);
	\fill[black!50!blue!50] (-4,-4) circle (.065);
	\fill[black!50!blue!50] (-5,-4) circle (.065);
	\fill[black!50!blue!50] (-6,-4) circle (.065);
	\fill[black!50!blue!50] (0,-5) circle (.065);
	\fill[black!50!blue!50] (1,-5) circle (.065);
	\fill[black!50!blue!50] (2,-5) circle (.065);
	\fill[black!50!blue!50] (3,-5) circle (.065);
	\fill[black!50!blue!50] (4,-5) circle (.065);
	\fill[black!50!blue!50] (5,-5) circle (.065);
	\fill[black!50!blue!50] (6,-5) circle (.065);
	\fill[black!50!blue!50] (-1,-5) circle (.065);
	\fill[black!50!blue!50] (-2,-5) circle (.065);
	\fill[black!50!blue!50] (-3,-5) circle (.065);
	\fill[black!50!blue!50] (-4,-5) circle (.065);
	\fill[black!50!blue!50] (-5,-5) circle (.065);
	\fill[black!50!blue!50] (-6,-5) circle (.065);
	\fill[black!50!blue!50] (0,-6) circle (.065);
	\fill[black!50!blue!50] (1,-6) circle (.065);
	\fill[black!50!blue!50] (2,-6) circle (.065);
	\fill[black!50!blue!50] (3,-6) circle (.065);
	\fill[black!50!blue!50] (4,-6) circle (.065);
	\fill[black!50!blue!50] (5,-6) circle (.065);
	\fill[black!50!blue!50] (6,-6) circle (.065);
	\fill[black!50!blue!50] (-1,-6) circle (.065);
	\fill[black!50!blue!50] (-2,-6) circle (.065);
	\fill[black!50!blue!50] (-3,-6) circle (.065);
	\fill[black!50!blue!50] (-4,-6) circle (.065);
	\fill[black!50!blue!50] (-5,-6) circle (.065);
	\fill[black!50!blue!50] (-6,-6) circle (.065);
	\draw[black, ultra thick, <->] (6.5,-5) -- (0,-5) -- (0,-4) -- (1,-3) -- (1,-2) -- (2,-1) -- (2,0) -- (3,1) -- (3,2) -- (4,3) -- (4,4) -- (5,5) -- (6.5,5);
	\draw[black, ultra thick, <->] (-3.5,-6.5) -- (-1,-4) -- (-1,-3) -- (0,-2) -- (0,-1) -- (1,0) -- (1,1) -- (2,2) -- (2,3) -- (3,4) -- (3,6.5);
	\draw[black, ultra thick, <->] (-5.5,-6.5) -- (-2,-3) -- (-2,-2) -- (-1,-1) -- (-1,0) -- (0,1) -- (0,2) -- (1,3) -- (1,6.5);
	\draw[black, ultra thick, <->] (-6.5,-5.5) -- (-3,-2) -- (-3,-1) -- (-2,0) -- (-2,1) -- (-1,2) -- (-1,6.5);
	\draw[black, ultra thick, <->] (-6.5,-3.5) -- (-4,-1) -- (-4,0) -- (-3,1) -- (-3,2) -- (-3,6.5);
	\draw[black, ultra thick, <->] (-6.5,-1.5) -- (-5,0) -- (-5,6.5);
	\draw[black, ultra thick, ->] (1,-3) -- (6.5,-3);
	\draw[black, ultra thick, ->] (2,-1) -- (6.5,-1);
	\draw[black, ultra thick, ->] (3,1) -- (6.5,1);
	\draw[black, ultra thick, ->] (4,3) -- (6.5,3);
	\draw[black, ultra thick, ->] (5,5) -- (5,6.5);
	\draw[black, ultra thick, ->] (0,-5) -- (-1.5,-6.5);
	\draw[black, ultra thick] (-1,-4) -- (0,-4);
	\draw[black, ultra thick] (0,-2) -- (1,-2);
	\draw[black, ultra thick] (1,0) -- (2,0);
	\draw[black, ultra thick] (2,2) -- (3,2);
	\draw[black, ultra thick] (3,4) -- (4,4);
	\draw[black, ultra thick] (-2,-3) -- (-1,-3);
	\draw[black, ultra thick] (-1,-1) -- (0,-1);
	\draw[black, ultra thick] (0,1) -- (1,1);
	\draw[black, ultra thick] (1,3) -- (2,3);
	\draw[black, ultra thick] (-3,-2) -- (-2,-2);
	\draw[black, ultra thick] (-2,0) -- (-1,0);
	\draw[black, ultra thick] (-1,2) -- (0,2);
	\draw[black, ultra thick] (-4,-1) -- (-3,-1);
	\draw[black, ultra thick] (-3,1) -- (-2,1);
	\draw[black, ultra thick] (-5,0) -- (-4,0);
	\draw[black, ultra thick, dotted] (6.65,5) -- (7,5);
	\fill[black] (7.1,5) circle (.08);
	\draw[black, ultra thick, dotted] (6.65,3) -- (7,3);
	\fill[black] (7.1,3) circle (.08);
	\draw[black, ultra thick, dotted] (6.65,1) -- (7,1);
	\fill[black] (7.1,1) circle (.08);
	\draw[black, ultra thick, dotted] (6.65,-1) -- (7,-1);
	\fill[black] (7.1,-1) circle (.08);
	\draw[black, ultra thick, dotted] (6.65,-3) -- (7,-3);
	\fill[black] (7.1,-3) circle (.08);
	\draw[black, ultra thick, dotted] (6.65,-5) -- (7,-5);
	\fill[black] (7.1,-5) circle (.08);
	\draw[black, ultra thick, dotted] (5,6.65) -- (5,7);
	\fill[black] (5,7.1) circle (.08);
	\draw[black, ultra thick, dotted] (3,6.65) -- (3,7);
	\fill[black] (3,7.1) circle (.08);
	\draw[black, ultra thick, dotted] (1,6.65) -- (1,7);
	\fill[black] (1,7.1) circle (.08);
	\draw[black, ultra thick, dotted] (-1,6.65) -- (-1,7);
	\fill[black] (-1,7.1) circle (.08);
	\draw[black, ultra thick, dotted] (-3,6.65) -- (-3,7);
	\fill[black] (-3,7.1) circle (.08);
	\draw[black, ultra thick, dotted] (-5,6.65) -- (-5,7);
	\fill[black] (-5,7.1) circle (.08);
	\draw[black, ultra thick, dotted] (-1.6,-6.6) -- (-1.85,-6.85);
	\fill[black] (-1.925,-6.925) circle (.08);
	\draw[black, ultra thick, dotted] (-3.6,-6.6) -- (-3.85,-6.85);
	\fill[black] (-3.925,-6.925) circle (.08);
	\draw[black, ultra thick, dotted] (-5.6,-6.6) -- (-5.85,-6.85);
	\fill[black] (-5.925,-6.925) circle (.08);
	\draw[black, ultra thick, dotted] (-6.6,-1.6) -- (-6.85,-1.85);
	\fill[black] (-6.925,-1.925) circle (.08);
	\draw[black, ultra thick, dotted] (-6.6,-3.6) -- (-6.85,-3.85);
	\fill[black] (-6.925,-3.925) circle (.08);
	\draw[black, ultra thick, dotted] (-6.6,-5.6) -- (-6.85,-5.85);
	\fill[black] (-6.925,-5.925) circle (.08);
	\end{tikzpicture}
	};
	(32,32)*+{{\color{blue}{\!1}}};
	(-28,2)*+{{\color{blue}{\!2}}};
	(2,-28)*+{{\color{blue}{\!3}}};
	(12,22)*+{{\color{blue}{\!4}}};
	(-18,-8)*+{{\color{blue}{\!5}}};
	(12,-8)*+{{\color{blue}{\!6}}};
	(-8,12)*+{{\color{blue}{\!7}}};
	(-8,-18)*+{{\color{blue}{\!8}}};
	(22,12)*+{{\color{blue}{\!9}}};
	(2,2)*+{{\color{blue}{\!10}}};
	\end{xy}
	$}
	$$
\caption{Galois-equivariant tropicalization of the degree-$6$ curve ${X}\subset Y_{\Sigma}$ in Example \ref{example: Brauer-Severi C}. Note the $\ZZ/3\ZZ$-symmetry under the action in Figures \ref{figure: example of action} and \ref{figure: example of tropical action}.}
\label{figure: example of tropicalization}
\end{figure}
	In the language of \cite[Definitions 7, 10, \& 11]{IKMZ:16}, $\text{Trop}(X,\imath)$ is a smooth regular projective $\QQ$-tropical variety inside $\text{Trop}(\PP^{2}_{\!\!\sm{\CC(\!(t^{\QQ})\!)}})\cong\TT\PP^{2}$. The dimensions of its nonzero tropical cellular cohomology groups $H^{p,q}_{\text{trop}}\big(\text{Trop}(X,\imath),\QQ\big)$ sit in a symmetric Hodge diamond
	$$
	\begin{array}{ccc}
	&
	\!1\!
	\\
	10\!\!
	&&
	\!\!\!10
	\\
	&
	\!1.\!\!
	\end{array}
	$$
The Galois group $G^{op}$ acts trivially on the degree-$(0,0)$ and -$(1,1)$ homology groups. Its action on the degree-$(1,0)$ and -$(0,1)$ groups factors through a $\ZZ/3\ZZ$-action that permutes the elements of a basis in each case. In degree-$(0,1)$, these basis elements can be taken as the classes of the $1$-cycles labeled ``1" through ``10" in Figure \ref{figure: example of tropicalization}, say with counter-clockwise orientations. With respect to this basis, the $G^{\op}$-action on $H_{0,1}^{\text{trop}}\big(\text{Trop}(X,\imath),\QQ\big)$ factors through the $\ZZ/3\ZZ$-action generated by the action of the $10\times 10$-matrix
	$$
	\begin{xy}
	(0,0)*+{
	\left(
	\begin{array}{c|c|c|c}
	\!\!\!\!
	\begin{xy}
	(0,0)*+{
	\begin{smallmatrix}0 & 0 & 1 \\ 1 & 0 & 0 \\ 0 & 1 & 0\end{smallmatrix}
	};
	\end{xy}
	\!\!\!\!\!
	&
	\multicolumn{3}{c}{}
	\\\cline{0-1}
	\multicolumn{1}{c|}{}
	&
	\!\!\!\!
	\begin{xy}
	(0,0)*+{
	\begin{smallmatrix}0 & 0 & 1 \\ 1 & 0 & 0 \\ 0 & 1 & 0\end{smallmatrix}
	};
	\end{xy}
	\!\!\!\!\!
	&
	\multicolumn{2}{c}{}
	\\\cline{2-3}
	\multicolumn{2}{c|}{}
	&
	\!\!\!\!
	\begin{xy}
	(0,0)*+{
	\begin{smallmatrix}0 & 0 & 1 \\ 1 & 0 & 0 \\ 0 & 1 & 0\end{smallmatrix}
	};
	\end{xy}
	\!\!\!\!\!
	&
	\\\cline{3-4}
	\multicolumn{3}{c|}{}
	&
	\!\!\!\!
	\begin{xy}
	(0,0)*+{
	\begin{smallmatrix}1\end{smallmatrix}
	};
	\end{xy}
	\!\!\!\!\!
	\end{array}
	\right)
	};
	(10,7)*+{\mbox{{\smaller\smaller $0$}}};
	(-7,-9)*+{\mbox{{\smaller\smaller $0$}}};
	\end{xy}.
	$$
	
	Under the isomorphism provided by \cite[Theorem 1]{IKMZ:16}, the $\ZZ/3\ZZ$-actions on tropical cohomology groups' $\QQ$-linear duals $H_{0,1}^{\text{trop}}\big(\text{Trop}(X,\imath),\QQ\big)^{\ast}$ and $H_{1,0}^{\text{trop}}\big(\text{Trop}(X,\imath),\QQ\big)^{\ast}$ induce $G$-representations on the respective weight-$0$ and -$1$ associated graded components of the weight filtration for the limit mixed Hodge structure of the nearby fiber associated to the family of projective varieties cut out by \eqref{equation: genus 3 Laurent polynomial} over the punctured unit disk inside $\text{Spec}_{\ \!}\CC[t]$.
\end{example}
	
\end{subsection}

\end{section}


\vskip .5cm
\begin{section}{An extensive supply of Galois-equivariant toric embeddings}\label{section: Cox rings, characteristic spaces, and Galois-equivariant embeddings}
	
	We explain how to construct an extensive supply of Galois-equivariant closed embeddings. Throughout \S\ref{section: Cox rings, characteristic spaces, and Galois-equivariant embeddings}, let $X$ be a variety over $K=K^{\text{sep}}_{0}$, without any assumption that $K_{0}$ is non-Archimedean.


\begin{subsection}{Algorithms for constructing non-equivariant toric embeddings}\label{subsection: algorithms for constructing non-equivariant toric embeddings}
	In order to prove Theorems \ref{theorem: main B} and \ref{theorem: main C}, we need to construct Galois-equivariant toric embeddings $X\mono Y_{\Sigma}$ that realize pre-specified rational functions on ${X}$ as pullbacks of characters on the dense torus in $Y_{\Sigma}$. To this end, we provide three constructions of toric embeddings in the non-equivariant setting: an embedding for projective varieties, J. W\l odarczyk's algorithm, and the algorithm of F. Berchtold and J. Hausen. W\l odarczyk's algorithm is sufficient when $K_{0}$ is a perfect field, but Berchtold and Hausen's algorithm, which lacks some of the generality of W\l odarczyk's, provides a far more explicit picture of what is happening and provides the user with tools that will be crucial (we believe) to applications.
	
\begin{ssubsection}\label{ssubsection: Payne's algorithm}
{\bf Payne's lemma for projective varieties.}
\normalfont
	In \cite{P09}, S. Payne shows how to construct systems of closed embeddings of a given quasiprojective variety into quasiprojective toric varieties such that the inverse limit of all tropicalizations in these systems is homeomorphic to the Berkovich analytification. A key lemma in the course of his proof is the following:
\end{ssubsection}

\begin{lemma}[{\it \cite[Lemma 4.3]{P09}}]\label{theorem: P embedding}
\normalfont
	Assume that ${X}$ is a projective ${K}$-variety. Given an ample effective divisor $E$ on ${X}$ with complement $U\subset{X}$, given a closed subvariety contained in the support of $E$, and given any set $f_{1},\dots,f_{r}$ of generators of the coordinate ring of $U$, there exists a closed toric embedding
	$$
	\imath:{X}\mono \PP^{n}_{\!\!\sm{{K}}}\ \Def\ \text{Proj}{K}[y_{0},\dots,y_{n}]
	$$
such that $U=\imath^{-1}\big(D(y_{0})\big)$, such that $V$ is the preimage of some coordinate linear subspace in $\PP^{n}_{\!\!\sm{{K}}}$, and such that each $f_{i}$ is the pullback $f_{i}=\imath^{\#}(y_{i})$.
\hfill
$\square$
\end{lemma}
	
\begin{ssubsection}\label{ssubsection: Wlodarczyk's algorithm}
{\bf W\l odarczyk's algorithm for A\textsubscript{2}-varieties.}
\normalfont
	Recall from Remark \ref{remark: A2} that an A\textsubscript{2}-variety is a variety in which every pair of points lies in some affine open subset of the variety.
	
	Let $K^{\text{alg}}\supset K$ denote an algebraic closure of $K$ (and thus an algebraic closure of $K_{0}$). Recall that $K_{0}$ is a {\em perfect} field if the inclusion ${K}=K^{\text{sep}}_{0}\subset K^{\text{alg}}_{0}=K^{\text{alg}}$ is an isomorphism. In this case $X_{K^{\text{alg}}}=X$. All fields are perfect in characteristic $0$.
	
	Suppose that ${X}_{K^{\text{alg}}}$ is a normal A\textsubscript{2}-variety over $K^{\text{alg}}$. In his seminal 1993 paper \cite{Wlod}, J. W\l odarczyk details an algorithm that takes as input any finite collection $\mathfrak{F}$ of nonzero rational functions on $X_{K^{\text{alg}}}$. The algorithm enlarges the collection $\mathfrak{F}$ through step-by-step adjustments of a corresponding set of divisors on $X_{K^{\text{alg}}}$. The final enlarged collection $\mathfrak{F}'$ satisfies properties that make it a generating set for the mutliplicative group of restrictions of characters under a closed embedding of $X_{K^{\text{alg}}}$ into a toric $K^{\text{alg}}$-variety $Y_{\Sigma,K^{\text{alg}}}$, and one can reconstruct $Y_{\Sigma,K^{\text{alg}}}$ from $\mathfrak{F}'$.
	
\begin{remark}
\normalfont
	W\l odarczyk works over the algebraic closure $K^{\text{alg}}$ because it is important in his algorithm that all classes in the quotient $K^{\text{alg}}(X)^{\times}\big/(K^{\text{alg}})^{\times}$ be torsion free \cite[bottom of p. 710]{Wlod}. If we take $K_{0}$ to be the non-perfect field $K_{0}=\FF_{\!p}(\!(t)\!)$ for instance, then on the $K=\FF_{\!p}(\!(t)\!)^{\text{sep}}$-variety $X=\text{Spec}_{\ \!}\FF_{\!q}(\!(t)\!)^{\text{sep}}[x]\big/(x^{p}-t)$ we have a non-trivial class $[x]\in \FF_{\!p}(\!(t)\!)^{\text{sep}}(X)^{\times}\big/\big(\FF_{\!p}(\!(t)\!)^{\text{sep}}\big)^{\times}$ with $[x]^{p}=[1]$.
\end{remark}
	
	Using W\l odarczyk's algorithm, P. Gross, S. Payne, and the present author show the following:
\end{ssubsection}

\begin{theorem}[{\it \cite[Theorem 4.2]{FGP}}]\label{theorem: FGP embedding}
\normalfont
	Assume that ${X}_{K^{\text{alg}}}$ is a normal A\textsubscript{2}-variety. Given affine open subvarieties $U_{1},\dots,U_{r}\subset X_{K^{\text{alg}}}$ and any finite collection $\mathfrak{F}=\mathfrak{F}_{1}\sqcup\cdots\sqcup\mathfrak{F}_{r}$ of rational functions such that each $f\in\mathfrak{F}_{i}$ is regular on $U_{i}$, there exists a closed embedding
	$$
	\imath:X_{K^{\text{alg}}}\mono Y_{\Sigma}
	$$
into a toric ${K^{\text{alg}}}$-variety $Y_{\Sigma,K^{\text{alg}}}$ such that for each $1\le i\le r$ we have $U_{i}=\imath^{-1}(Y_{\sigma_{i},K^{\text{alg}}})$ for some $\sigma_{i}\in\Sigma$, and such that each $f\in\mathfrak{F}_{i}$ is the pullback $\imath^{\#}(\chi^{u})$ of a character $\chi^{u}\in S_{\sigma_{i}}$.
\hfill
$\square$
\end{theorem}
	
\begin{ssubsection}\label{ssubsection: Berchtold and Hausen's algorithm}
{\bf Berchtold and Hausen's algorithm for ``A\textsubscript{2}-Mori dream spaces" in characteristic 0.}
\normalfont
	In the last decade and a half, F. Berchtold and J. Hausen \cite{BH:04} \cite{BH:07} \cite{Haus:08} have clarified large parts of the geometry behind W\l odarczyk's original algorithm using the language of Cox rings. See Hausen's book with I. Arzhantsev, U. Derenthal, and A. Laface \cite[\S1-\S3]{ADHL} for a detailed exposition. Berchtold and Hausen employ combinatorial objects called ``bunches of cones," dual to fans, inside the Cox ring of ${X}$. The clarity that Berchtold and Hausen's construction affords comes at a cost: one must assume that $\text{char}_{\ \!}K=0$ (thus that $K_{0}$ is perfect) and that ${X}$ satisfies an A\textsubscript{2} version of the Mori dream space hypotheses. Specifically, $X$ must be an irreducible factorial A${}_{2}$-variety over $K^{\text{alg}}={K}$ with $H^{0}({X},\mathscr{O}^{\times}_{\!\sm{{X}}})=K^{\times}$ such that the ideal class group and Cox ring of ${X}$ are finitely generated. These conditions hold, for instance, for any $G^{\op}$-twisted proper (not-necessarily projective) toric ${K}$-variety in characteristic $0$. Using a finite set of generators $\mathfrak{F}$ of ${X}$'s Cox ring $R_{{X}}$, Berchtold and Hausen construct a characteristic space over ${X}$ that comes equipped with a canonical closed embedding into a toric ${K}$-variety $Y_{\widehat{\Sigma}}$. Taking a quotient of both $Y_{\widehat{\Sigma}}$ and the characteristic space lying inside it, one obtains a closed embedding
	$
	\imath:{X}\mono Y_{\Sigma}
	$
into a toric ${K}$-variety $Y_{\Sigma}$.
	
	In Theorem \ref{theorem: BH embedding} below, we prove an analogue of Theorem \ref{theorem: FGP embedding} that uses Berchtold and Hausen's algorithm in place of W\l odarczyk's algorithm. The proof makes heavy use of notation and terminology from \cite{ADHL}. It is the only place in the present text where we make use of this notation, and readers interested only in the proofs of Theorems \ref{theorem: main A}, \ref{theorem: main B}, and \ref{theorem: main C} may want to read the statement of Theorem \ref{theorem: BH embedding} and then skip ahead to \S\ref{subsection: from non-equivariant to equivariant toric embeddings}.
\end{ssubsection}

\begin{theorem}\label{theorem: BH embedding}
\normalfont
	Assume that $\text{char}_{\ \!}K_{0}=0$, that ${X}$ is an A${}_{2}$-variety over ${K}$ with $H^{0}({X},\mathscr{O}^{\times}_{\!\sm{{X}}})=K^{\times}$, and that the ideal class group $\text{Cl}(X)$ and Cox ring $R_{X}$ of $X$ are finitely generated. Then there exists a finite affine open cover $\mathscr{U}=\{U_{1},\dots,U_{r}\}$ of ${X}$, consisting of images of standard affine opens in $\text{Spec}_{\ \!}R_{X}$, and satisfying the following property: For any function $f$ regular on some $U_{i}$, we can construct a closed embedding
	$$
	\imath':{X}\mono Y_{\Sigma'}
	$$
into a toric ${K}$-variety $Y_{\Sigma'}$ such that for each $1\le i\le r$, the affine open $U_{i}$ is the inverse image of a torus invariant open in $Y_{\Sigma'}$, and such that some $K^{\times}$-multiple of $f$ is the pullback of a character on the dense torus in $Y_{\Sigma'}$.
\end{theorem}
\begin{proof}
	By \cite[Theorem 1.5.3.7]{ADHL}, we can choose an ordered generating set $\mathfrak{F}=(f_{1},\dots,f_{r})$ of the Cox ring $R_{{X}}$ satisfying the hypotheses of \cite[\S3.2.1]{ADHL}. By \cite[Theorem 3.2.1.9.(ii)]{ADHL} and the A\textsubscript{2} hypothesis, we can find a maximal $\mathfrak{F}$-bunch $\Phi$. Via \cite[Construction 3.2.1.3]{ADHL}, this provides us with a characteristic space over ${X}$ and an affine open cover $\mathscr{U}=\{U_{\gamma_{I}}\}_{I\in\mathscr{I}}$ of ${X}$ indexed by the set $\mathscr{I}$ of relevant faces of $\Phi$. Using \cite[Construction 3.2.5.3 \& Proposition 3.2.5.4]{ADHL}, we obtain a closed embedding into a toric variety
	$$
	\imath_{\Phi}:{X}\mono Y_{\Sigma}
	$$
such that each $U_{\gamma_{I}}$ is of the form $\imath^{-1}(Y_{\sigma})$ for some $\sigma\in\Sigma$.
	
	We want to check that for each function $f$ regular on some $U_{\gamma_{I}}$, there exists another closed toric embedding $\imath':{X}\mono Y_{\Sigma'}$ such that each $U_{\gamma_{I}}$ is the inverse image of a torus-invariant affine open, and such that some $K^{\times}$-multiple of $f$ is the pullback of a character. Fix one such $U_{\gamma_{I}}$ and fix a regular function $f$ on $U_{\gamma_{I}}$. By construction, we can interpret $f$ as a degree-$0$ element
	$$
	f
	\ =\ 
	\frac{f_{0}}{(f^{u_{1}}_{1}\cdots f^{u_{r}}_{r})^{n}}
	\ \ \in\ \ 
	R_{{X}}\!\!\left[\frac{1}{f^{u_{1}}_{1}\cdots f^{u_{r}}_{r}}\right]
	$$
for some $n\ge0$ and for any choice of $u=(u_{1},\dots,u_{r})\in\gamma^{\circ}_{I}$, where $f_{0}$ is homogeneous of degree $n\ \!Q(u)\in\text{Cl}({X})$. By \cite[Theorem 1.5.3.7]{ADHL}, we can assume that $f_{0}$ is $\text{Cl}({X})$-prime.

	If $f_{0}$ is associated to one of the functions in $\mathfrak{F}$, then since $H^{0}({X},\mathscr{O}^{\times}_{\!\sm{{X}}})=K^{\times}$, after multiplying a scalar in $K^{\times}$ we have that $f_{0}\in\mathfrak{F}$. In this case, the fact that $f$ is degree-$0$ implies that it is the pullback under $\imath_{\Phi}$ of a character that is regular on the affine toric subvariety $Y_{P(\gamma^{\ast}_{I})}$.
	
	If $f_{0}$ is {\em not} associated to any function in $\mathfrak{F}$, then by the previous arguments, the ordered set of homogeneous elements
	$
	\mathfrak{F}'
	\ =\ 
	(f_{0},f_{1},\dots,f_{r})
	$
admits an $\mathfrak{F}'$-bunch $\Phi'$ that provides us with a new closed toric embedding $\imath_{\Phi''}:{X}\mono Y_{\Sigma''}$ for which $f$ becomes the pullback of a monomial. Then the product toric embedding
	$$
	\imath'\Def\imath_{\Phi}\times\imath_{\Phi''}:{X}\mono Y_{\Sigma\times\Sigma''}
	$$
satisfies the two requisite properties.
\end{proof}

\end{subsection}


\begin{subsection}{From non-equivariant to equivariant toric embeddings}\label{subsection: from non-equivariant to equivariant toric embeddings}
	
	Assume now and for the remainder of the text that $X$ arises as $X=X_{0,{K}}$ for some $K_{0}$-variety $X_{0}$ and $K=K^{\text{sep}}_{0}$. Equip $X$ with its canonical $G^{\op}$-action. Since every toric variety has a canonical $\ZZ$-model, the toric ${K}$-varieties produced by Lemma \ref{theorem: P embedding} when $X$ is projective, by Theorem \ref{ssubsection: Wlodarczyk's algorithm} when $K_{0}$ is perfect, and by Theorem \ref{ssubsection: Berchtold and Hausen's algorithm} when $\text{char}_{\ \!}K_{0}=0$ come at least with a canonical $G^{\op}$-actions on $X$ and $Y_{\Sigma}$. However, because each algorithm constructs closed embeddings over ${K}$, not $K_{0}$, the embeddings may not be $G$-equivariant with respect to the canonical $G^{\op}$-action. We get $G$-equivariant toric embeddings by running a second construction reminiscent of the construction of an induced representation.

\begin{theorem}\label{equation: making a toric embedding G-equivariant}
\normalfont
	For each closed embedding $\imath:{X}\mono Y_{\Sigma}$ into a toric ${K}$-variety $Y_{\Sigma}$ associated to a fan $\Sigma\subset N_{\RR}$ without $G^{\op}$-action, there exists a $G^{\op}$-equivariant closed embedding $\jmath:X\mono Y_{\Sigma'}$ into a $G^{\op}$-twisted toric variety $Y_{\Sigma'}$ and a (non-$G$-equivariant) surjective morphism of toric ${K}$ varieties $Y_{\Sigma'}\epi Y_{\Sigma}$ that fits into a commutative diagram of $K$-varieties
	\begin{equation}\label{equation: crucial projection}
	\begin{xy}
	(0,0)*+{\ {X}\ }="1";
	(17,8)*+{Y_{\Sigma'}}="2";
	(17,-8)*+{Y_{\Sigma}.\!}="3";
	{\ar@{^{(}->}^{\jmath} "1"; "2"};
	{\ar@{_{(}->}_{\imath} "1"; "3"};
	{\ar@{->>} "2"; "3"};
	\end{xy}
	\end{equation}
\end{theorem}
\begin{proof}
	The toric ${K}$-variety $Y_{\Sigma}$ is of the form ${K}\otimes_{K_{0}}Y_{\Sigma,K_{0}}$, where $Y_{\Sigma,K_{0}}$ denotes the toric $K_{0}$-variety associated to $\Sigma$. Equip $Y_{\Sigma}$ with its resulting canonical $G^{\op}$-action. Let $M$ denote the lattice of characters on the dense torus in $Y_{\Sigma}$. Define $\mathfrak{F}$ to be the multiplicative group
	$$
	\mathfrak{F}
	\ \Def\ 
	\big\{
	\imath^{\#}(\chi^{u}):u\in M
	\big\},
	$$
consisting of rational functions on ${X}$. The canonical $G^{\op}$-action $G^{\op}\ \lact\ {X}$ induces a $G$-action on the multiplicative semigroup ${K}({X})$ of rational functions. Because $M$ is a finitely generated free group, $\mathfrak{F}$ is a finitely generated multiplicative subsemigroup of $K(X)$. Hence the $G$-orbit $G\mathfrak{F}$ is itself a finitely generated multiplicative subsemigroup of ${K}({X})$. Consequently, there exists an open normal subgroup $Z\lhd G$ fixing $\mathfrak{F}$, i.e., such that $Zf=f$ for all $f\in\mathfrak{F}$. If we define
	$$
	H
	\ \Def\ 
	G/Z,
	$$
then the natural $G$-action on the orbit $G\mathfrak{F}$ factors through an $H$-action $H\ \lact\ G\mathfrak{F}$. Let $L\Def{K}^{Z}$ denote the fixed field of $Z$, so that $H=\text{Gal}(L/K_{0})$.
	
	If $Y_{\Sigma,L}$ denotes the toric $L$-variety associated to the fan $\Sigma$, then our construction of $L$ implies that the closed embedding $\imath:{X}\mono Y_{\Sigma}$ over $K$ is actually the pullback of a closed embedding
	$$
	\imath_{L}:X_{0,L}\mono Y_{\Sigma,L}.
	$$
Equip $X_{0,L}$ and $Y_{\Sigma,L}$ with their canonical $H^{\op}$-actions. For each $h\in H$, let $\imath^{h}_{L}:X_{0,L}\mono Y_{\Sigma,L}$ denote the pullback morphism fitting into the Cartesian diagram
	\begin{equation}\label{equation: pulled back inclusion}
	\begin{aligned}
	\begin{xy}
	(0,0)*+{X_{0,L}}="1";
	(20,0)*+{Y_{\Sigma,L}}="2";
	(0,-18)*+{X_{0,L}}="3";
	(20,-18)*+{Y_{\Sigma,L}.\!}="4";
	{\ar@{^{(}->}^{\imath^{h}_{L}\ } "1"; "2"};
	{\ar@{^{(}->}_{\imath_{L}} "3"; "4"};
	{\ar_{(g^{-1})^{\ast}} "1"; "3"};
	{\ar^{(g^{-1})^{\ast}} "2"; "4"};
	(5,-5)*+{\mbox{{\larger\larger\larger\larger $\lrcorner$}}};
	\end{xy}
	\end{aligned}
	\end{equation}
Define $N^{H}\Def\prod_{H}N$ to be $\#H$-fold Cartesian power of $N$ with factors indexed by $H$. Likewise, define $\Sigma^{H}\Def\prod_{H}\Sigma$ to be the fan in $N^{H}_{\RR}$ obtained as the $\#H$-fold Cartesian power of the fan $\Sigma$ with factors indexed by $H$. Let $Y_{\Sigma^{H}\!\!,\ \!L}$ be the associated toric $L$-variety. Then we have a closed embedding
	\begin{equation}\label{equation: newly constructed equivariant embedding}
	\jmath_{L}\Def\mbox{$\prod_{H}$}\ \imath^{h}_{L}\ :\ \ X_{0,L}\mono Y_{\Sigma^{H}\!\!,\ \!L}.
	\end{equation}
The cocharacter lattice $N^{H}$ comes with a natural $H^{\op}$-action that permutes the factors of $N^{H}$, and the fan $\Sigma^{H}$ is invariant under the induced $H^{\op}$-action $H^{\op}\ \lact\ N^{H}_{\RR}$. This $H^{\op}$-action on $N^{H}$ induces a twisted $H^{\op}$-action $H^{\op}\ \lact\ Y_{\Sigma^{H}\!\!,\ \!L}$,  which is to say that the diagram
	$$
	\begin{aligned}
	\begin{xy}
	(0,0)*+{Y_{\Sigma^{H}\!\!,\ \!L}}="1";
	(22,0)*+{Y_{\Sigma^{H}\!\!,\ \!L}}="2";
	(0,-15)*+{\text{Spec}_{\ \!}L}="3";
	(22,-15)*+{\text{Spec}_{\ \!}L}="4";
	{\ar^{h^{\!\ast}} "1"; "2"};
	{\ar_{h^{\!\ast}} "3"; "4"};
	{\ar@{->>} "1"; "3"};
	{\ar@{->>} "2"; "4"};
	\end{xy}
	\end{aligned}
	$$
commutes for each $h\in H$. 
	
	We claim that the closed embedding \eqref{equation: newly constructed equivariant embedding} is $H$-equivariant, i.e., that each diagram
	\begin{equation}\label{equation: equivariance check diagram}
	\begin{aligned}
	\begin{xy}
	(0,0)*+{X_{0,L}}="1";
	(20,0)*+{Y_{\Sigma^{H}\!\!,\ \!L}}="2";
	(0,-18)*+{X_{0,L}}="3";
	(20,-18)*+{Y_{\Sigma^{H}\!\!,\ \!L}}="4";
	{\ar@{^{(}->}^{\jmath_{L}\ } "1"; "2"};
	{\ar@{^{(}->}_{\jmath_{L}\ } "3"; "4"};
	{\ar_{h^{\ast}} "1"; "3"};
	{\ar^{h^{\ast}} "2"; "4"};
	\end{xy}
	\end{aligned}
	\end{equation}
commutes. Since $Y_{\Sigma^{H}\!\!,\ \!L}=\prod_{H}Y_{\Sigma,L}$, we can check this by verifying the universal property of each of the two composites in \eqref{equation: equivariance check diagram}. To this end, observe that for each $\ell\in H$, we have commutative diagrams
	$$
	\xymatrix{
	X_{0,L}
	\ar@{^{(}->}[r]_{\jmath_{L}\ \ }
	\ar@{^{(}->}@/^43pt/[rrr]^{h^{\!\ast}\imath^{(\ell h^{-1})}_{L}}
	&
	Y_{\Sigma^{H}\!\!,\ \!L}
	\ar[r]^{\mbox{$\sim$}}_{h^{\!\ast}}
	\ar@{->>}@/^18pt/[rr]^{h^{\!\ast}\text{pr}_{(\ell h^{-1})}\ \ \ \ \ }
	&
	Y_{\Sigma^{H}\!\!,\ \!L}
	\ar@{->>}[r]_{\text{pr}_{\ell}}
	&
	Y_{\Sigma,L}
	}
	$$
\begin{center}
and
\end{center} 
	$$
	\xymatrix{
	X_{0,L}
	\ar@{^{(}->}[r]^{h^{\!\ast}}_{\mbox{$\sim$}}
	&
	X_{0,L}
	\ar@{^{(}->}[r]_{\mbox{$\sim$}}^{\jmath_{L}}
	\ar@{_{(}->}@/_20pt/[rr]_{\imath^{\ell}_{L}}
	&
	Y_{\Sigma^{H}\!\!,\ \!L}
	\ar@{->>}[r]^{\text{pr}_{\ell}}
	&
	Y_{\Sigma,L}.
	}
	$$
Because the diagram \eqref{equation: pulled back inclusion} is Cartesian, the two composites are equal.

	Finally, define $Y_{\Sigma^{H}}$ to be the pullback of $Y_{\Sigma^{H}\!\!,\ \!L}$ to ${K}$, equipped with the resulting twisted $G^{\op}$-action, and define $\jmath:{X}\mono Y_{\Sigma^{H}}$ to be the pullback of $\jmath_{L}$ to ${K}$. Then $G$-equivariance of $\jmath$ follows immediately from $H$-equivariance of $\jmath_{L}$.
\end{proof}
	
\begin{corollary}\label{corollary: projection equivariant embeddings}
\normalfont
	Assume that ${X}$ is projective. Then for each integer $n>0$ and each non-equivariant closed embedding $\imath:{X}\mono\PP^{n}_{\!\!\sm{K}}$, where $\PP^{n}_{\!\!\sm{K}}$ is equipped with its standard toric structure, there exists a $G$-equivariant closed embedding
	$$
	\jmath':{X}\mono\PP^{m}_{\!\!\sm{{K}}}
	$$
that factors as
	$$
	{X}\xymatrix{{}\ar@{^{(}->}[r]^{\jmath}&{}}\prod_{H}\PP^{n}_{\!\!\sm{{K}}}\xymatrix{{}\ar@{^{(}->}[r]^{\text{Seg}}&{}}\PP^{m}_{\!\!\sm{{K}}}
	$$
where $\jmath$ is a $G$-equivariant closed embedding constructed from $\imath$ as in the proof of Proposition \ref{equation: making a toric embedding G-equivariant}, where $\text{Seg}$ is the $\#H$-fold Segre embedding, and where $\PP^{m}_{\!\!\sm{{K}}}$ has the $G^{\op}$-twisted action induced by the $H$-action on the factors of $\prod_{H}\PP^{m}_{\!\!\sm{{K}}}$.
\hfill
$\square$
\end{corollary}

\begin{remark}\label{remark: large supply of G-representations}
{\bf A large supply of {\em G}-representations.}
	Assume that ${X}$ is projective. Then by Lemma \ref{theorem: P embedding}, we can construct a large system of closed toric embeddings $\imath:{X}\mono\PP^{n}_{\!\!\sm{{K}}}$. Each embedding realizes a given set of functions on the complement of an effective ample divisor on ${X}$ as the pullbacks of coordinate linear functions on $\PP^{n}_{\!\!\sm{{K}}}$. By Corollary \ref{corollary: projection equivariant embeddings}, we can turn any given choice of $\imath$ into a $G$-equivariant closed embedding $\jmath':{X}\mono\PP^{m}_{\!\!\sm{{K}}}$. Thus we have an extensive supply of $G$-equivariant closed embeddings of ${X}$ into $G^{\op}$-twisted toric projective spaces. If we assume that $K_{0}$ is Henselian, then Theorem \ref{theorem: main D} turns this into a extensive supply of $G$-representations in tropical cellular cohomology groups.
\end{remark}

\end{subsection}


\end{section}


\vskip .5cm
\begin{section}{Proofs of the Main Theorems}
\label{section: proofs of the main theorems}


\begin{subsection}{Criterion for a {\em G}-equivariant homeomorphism: Theorem \ref{theorem: main A}}\label{subsection: proof of Theorem A}
	We begin by proving Theorem \ref{theorem: main A}, which says that Conditions ($\Pi$) and ($\star$) together provide a criterion for checking if a system of $G$-equivariant toric varieties has enough closed embeddings to give an affirmative answer to Question \ref{questions: recovering Berkovich analytification tropically}. The proof amounts to an appeal, through cofinality, to the results of \cite{FGP}.

\begin{proof}[{\it Proof of Theorem \ref{theorem: main A}}]
	Let $\mathcal{S}_{G}$ be a system of $G$-equivariant toric embeddings of ${X}$ satisfying Conditions ($\Pi$) and ($\star$) in Definition \ref{defintion: G-equivariant system}. There is a forgetful functor
	\begin{equation}\label{equation: forgetful functor}
	\mathcal{S}_{G}\longrightarrow\bold{TEmb}({X})
	\end{equation}
that takes each $G$-equivariant closed toric embedding in $\mathcal{S}_{G}$ to its underlying closed toric embedding. Note that this forgetful functor \eqref{equation: forgetful functor} will not be fully faithful, since a given toric ${K}$-variety $Y_{\Sigma}$ can have several distinct twisted $G^{\op}$-actions. Let $\mathcal{S}$ denote the system of toric embeddings of ${X}$ given by the image of the forgetful functor \eqref{equation: forgetful functor}. The fact that $\mathcal{S}_{G}$ satisfies Conditions ($\Pi$) and ($\star$) of Definition \ref{defintion: G-equivariant system} implies that $\mathcal{S}$ satisfies the hypotheses of \cite[Theorem 1.1]{FGP}. Thus the map
	$$
	{X}^{\an}
	\longrightarrow
	\varprojlim_{\imath\in\mathcal{S}}\text{Trop}({X},\imath)
	$$
is a homeomorphism of topological spaces (without $G^{\op}$-actions). The functor
	$$
	\mathcal{S}_{G}\epi\mathcal{S}
	$$
is final (or what is also called ``co-cofinal" in \cite[Definition 2.5.1.(ii)]{KS:06}). By \cite[Proposition 2.5.2.(ii)]{KS:06}, this implies that the forgetful functor induces a homeomorphism of the underlying topological spaces
	$$
	\varprojlim_{\imath\in\mathcal{S}}\text{Trop}({X},\imath)
	\ \xrightarrow{\ \ \sim\ \ }\ 
	\varprojlim_{\jmath\in\mathcal{S}_{G}}\text{Trop}_{G}({X},\jmath).
	$$
Hence the map \eqref{equation: canonical G-equivariant map} is a homeomorphism. Finally, the map \eqref{equation: canonical G-equivariant map} is $G$-equivariant by Corollary \ref{corollary: equivariance of tropicalization map}. 
\end{proof}

\end{subsection}


\begin{subsection}{Proofs of Theorems \ref{theorem: main B} and \ref{theorem: main C}}\label{subsection: proof of A, B, and C}
	 Theorem \ref{theorem: main A} provides us with a criterion for verifying whether or not a given system $\mathcal{S}_{G}$ of $G$-equivariant toric embeddings induces a $G$-equivariant homeomorphism \eqref{equation: canonical G-equivariant map}. In order to construct $G$-equivariant systems of toric embeddings satisfying this criterion, we employ any one of the algorithms of \S\ref{subsection: algorithms for constructing non-equivariant toric embeddings}.

\vskip .3cm

\begin{proof}[{\it Proof of Theorem \ref{theorem: main B}}]
	Fix an affine open cover $\{U_{1},\dots,U_{r}\}$ provided by Lemma \ref{theorem: P embedding} if $X$ is quasiprojective, or as in Theorem \ref{theorem: FGP embedding} if $K_{0}$ is perfect and $X$ admits at least one closed embedding into a toric ${K}$-variety (using the fact that any ambient toric variety is normal and A\textsubscript{2}), or as in Theorem \ref{theorem: BH embedding} if $\text{char}_{\ \!}K_{0}=0$ and $X$ is an A${}_{2}$-variety over ${K}$ with $H^{0}({X},\mathscr{O}^{\times}_{\!\sm{{X}}})=K^{\times}$ and the ideal class group and Cox ring of ${X}$ are finitely generated. Let $\mathfrak{F}$ denote the set of all rational functions on ${X}$ regular on at least one of the affine opens $U_{i}$. For each $f\in \mathfrak{F}$ and each choice of open $U_{i}$ with $f\in {K}[U_{i}]$, construct a non-equivariant toric embedding $\imath_{f}:{X}\mono Y_{\Sigma_{f}}$ as in Lemma \ref{theorem: P embedding}, Theorem \ref{theorem: FGP embedding}, or Theorem \ref{theorem: BH embedding}, respectively. Use $\imath_{f}$ and Theorem \ref{equation: making a toric embedding G-equivariant} to build a $G$-equivariant closed toric embedding
	\begin{equation}\label{equation: ultimate thing associated to f}
	\jmath_{f}:{X}\mono Y_{\Sigma^{H}_{f}}.
	\end{equation}
Note that the existence of the projection $Y_{\Sigma^{H}_{f}}\epi Y_{\Sigma}$ fitting into an instance of the commutative diagram \eqref{equation: crucial projection} implies that $U_{i}$ is the inverse image, under $\jmath_{f}$, of a torus-invariant open subscheme in $Y_{\Sigma^{H}_{f}}$, and that $f$ is the inverse image $f=\jmath^{\#}_{f}(a\chi^{u})$ of a monomial on the dense torus in $Y_{\Sigma^{H}_{f}}$.
	
	Let $\bold{Sep}_{G}(X;\mathfrak{F})$ denote the discrete subcategory of $\bold{TEmb}_{G}({X})$ with objects consisting of all the $G$-equivariant closed embeddings \eqref{equation: ultimate thing associated to f} that we've constructed as $(f,U_{i})$ ranges over $f\in\mathfrak{F}$ and $U_{i}$ such that $f\in {K}[U_{i}]$. Then $\bold{Sep}_{G}(X;\mathfrak{F})$ satisfies Condition ($\star$) of Definition \ref{defintion: G-equivariant system}. Any system $\mathcal{S}_{G}$ of $G$-equivariant toric embeddings containing $\bold{Sep}_{G}(X;\mathfrak{F})$ and satisfying Condition ($\Pi$) of Definition \ref{defintion: G-equivariant system} satisfies the theorem, with $\bold{TEmb}_{G}({X})$ being an example of one such system.
\end{proof}

\vskip .3cm

\begin{proof}[{\it Proof of Theorem \ref{theorem: main C}}]
	Let $x$ be a point in $X^{\an}$ with finite $G^{\op}$-orbit $G^{\op}x$. For each pair of distinct points $y$ and $z$ in $G^{\op}x$, \cite[Proposition 3.2]{FGP} shows that one of the closed embeddings $\imath_{yz}:{X}\mono Y_{\Sigma_{yz}}$ provided by Lemma \ref{theorem: P embedding}, Theorem \ref{theorem: FGP embedding}, or Theorem \ref{theorem: BH embedding} satisfies
	$$
	\text{trop}(x')\ \ne\ \text{trop}(x'')
	\ \ \ \ \mbox{inside}\ \ \ \ 
	\text{Trop}({X},\imath_{yz}).
	$$
Applying Theorem \ref{equation: making a toric embedding G-equivariant}, we obtain a $G$-equivariant toric embedding
	$$
	\jmath_{yz}:{X}\mono Y_{\Sigma_{yz}}
	$$
such that $\text{trop}(x')\ne\text{trop}(x'')$ inside $\text{Trop}_{G}({X},\jmath_{yz})$. Let $\text{S}^{2}_{\ \!}G^{\op}x$ denote the set of all unordered pairs of distinct points in $G^{\op}x$. Then we can form a new $G$-equivariant closed toric embedding as the product
	$$
	\jmath\Def\underset{\text{S}^{2}_{\ \!}G^{\textsf{op}}x}{\prod_{\{y,z\}\ \text{in}}}\jmath_{yz}
	\ :\ \ 
	{X}
	\mono
	\underset{\text{S}^{2}_{\ \!}G^{\textsf{op}}x}{\prod_{\{y,z\}\ \text{in}}}Y_{\Sigma_{yz}}.
	$$
Each projection
	$$
	\text{pr}_{yz}:\underset{\text{S}^{2}_{\ \!}G^{\textsf{op}}x}{\prod_{\{y,z\}\ \text{in}}}Y_{\Sigma_{yz}}\epi Y_{\Sigma_{yz}}
	$$
induces a morphism of $G$-equivariant tropicalizations $\text{Trop}_{G}(\prod_{\{y,z\}}Y_{\Sigma_{yz}})\longrightarrow\text{Trop}_{G}(Y_{\Sigma_{yz}})$ that restricts to a morphism of $G$-equivariant tropicalizations $\text{Trop}_{G}({X},\jmath)\longrightarrow\text{Trop}_{G}({X},\jmath_{yz})$. Because each pair of distinct points $y,z\in G^{\op}x$ is separated in its corresponding tropicalization $\text{Trop}({X},\jmath_{yz})$, we conclude that the $G$-equivariant composite
	$$
	G^{\op}x
	\mono
	{X}^{\an}
	\epi
	\text{Trop}({X},\jmath)
	$$
is a bijection onto its image.
\end{proof}

\end{subsection}

\end{section}


\vskip 1cm

\bibliography{All_Galois_action}
\bibliographystyle{plain}

\end{document}